\RequirePackage{fix-cm}
\documentclass[smallextended]{svjour3} 

\usepackage{graphicx}
\usepackage{color,verbatim}
\usepackage{epstopdf}

\usepackage{subfigure}
\usepackage{amsmath}
\usepackage{amsfonts,amssymb}
\usepackage{hyperref}
\usepackage{array}  
\usepackage{algorithm}
\usepackage{algorithmic}
\newtheorem{thm}{Theorem}[section]
\newtheorem{pro}[thm]{Proposition}
\newtheorem{lem}[thm]{Lemma}
\newtheorem{cor}[thm]{Corollary}
\newtheorem{ass}[thm]{Assumption}
\newtheorem{defi}[thm]{Definition}

\newcommand{\nn}{\nonumber}
\newcommand{\be}{\begin{equation}}
\newcommand{\ee}{\end{equation}}
\newcommand{\ba}{\begin{array}}
	\newcommand{\ea}{\end{array}}
\newcommand{\bea}{\begin{eqnarray}}
\newcommand{\eea}{\end{eqnarray}}
\newcommand{\bal}{\begin{alg}}
	\newcommand{\eal}{\end{alg}}
\newcommand{\ble}{\begin{lem}}
	\newcommand{\ele}{\end{lem}}
\newcommand{\bco}{\begin{cor}}
	\newcommand{\eco}{\end{cor}}
\newcommand{\bde}{\begin{defi}}
	\newcommand{\ede}{\end{defi}}
\newcommand{\bpr}{\begin{pro}}
	\newcommand{\epr}{\end{pro}}
\newcommand{\bas}{\begin{ass}}
	\newcommand{\eas}{\end{ass}}
\newcommand{\bex}{\begin{example}}
	\newcommand{\eex}{\end{example}}
\newcommand{\bfig}{\begin{figure}}
	\newcommand{\efig}{\end{figure}}
\newcommand{\reff}[1]{(\ref{#1})}
\newcommand{\refa}[1]{Assumption\ \ref{#1}}
\newcommand{\refl}[1]{Lemma\ \ref{#1}}
\newcommand{\reft}[1]{Theorem\ \ref{#1}}

\newcommand{\refal}[1]{Algorithm\ \ref{#1}}

\def\alglist{
	\begin{list}{Step 1}
		{\setlength{\leftmargin}{0.5 in}\setlength{\labelwidth}{0.7 in}}
	}
	\def\eli{\end{list}}

\def\na{\nabla}

\def\De{\Delta}

\def\la{\lambda}

\def\hf{\frac{1}{2}}

\def\alp{\alpha}

\def\st{\hbox{s.t.}}
\def\diag{\hbox{diag}\,}

\begin{document}

\title{A primal-dual interior-point relaxation method with
	global and rapidly local convergence for nonlinear programs}

\titlerunning{A primal-dual interior-point relaxation method...y} 

\author{Xin-Wei Liu         \and
	Yu-Hong Dai
	 \and
	Ya-Kui Huang 
}

\institute{Xin-Wei Liu \at
	School of Sciences, Hebei University of Technology, Tianjin 300401, China \\
	\email{mathlxw@hebut.edu.cn}           
	\and	
	Yu-Hong Dai \at
	LSEC, ICMSEC, Academy of Mathematics and Systems Science, Chinese Academy of Sciences, Zhongguancun East Road No. 55, Beijing 100190, China\\
	\email{dyh@lsec.cc.ac.cn}   	
	\and
	Ya-Kui Huang \at
	Institute of Mathematics, Hebei University of Technology, Tianjin 300401, China\\
	\email{hyk@hebut.edu.cn}
}

\date{Received: date / Accepted: date}

\maketitle


\begin{abstract}
	Based on solving an equivalent parametric equality constrained mini-max problem of the classic logarithmic-barrier subproblem, we present a novel primal-dual interior-point relaxation method for nonlinear programs with general equality and nonnegative constraints. In each iteration, our method approximately solves the KKT system of a parametric equality constrained mini-max subproblem, which avoids the requirement that any primal or dual iterate is an interior-point. The method	has some similarities to the warmstarting interior-point methods in relaxing the interior-point requirement and is easily extended for solving problems with general inequality constraints. In particular, it has the potential to circumvent the jamming difficulty that appears with many interior-point methods for nonlinear programs and improve the ill conditioning of existing primal-dual interior-point methods as the barrier parameter is small. A new smoothing approach is introduced to develop our relaxation method and promote convergence of the method. Under suitable conditions, it is proved that our method can be globally convergent and locally quadratically convergent to the KKT point of the original problem. The preliminary numerical results on a well-posed problem for which many interior-point methods fail to find the minimizer and a set of test problems from the CUTEr collection show that our method is efficient.

\end{abstract}

\section{Introduction}\label{sec1}

We consider the nonlinear programs with the form
\begin{align}
\hbox{minimize}\quad(\min) &\quad f(x) \label{prob1-1}\\
\hbox{subject to}\quad(\st)   &\quad h(x)=0, \quad x\ge 0, \label{prob1-2}
\end{align}
where $x\in\Re^n$, $f: \Re^n\rightarrow \Re$ and $h: \Re^n\rightarrow \Re^m$ are twice continuously differentiable real-valued functions defined on $\Re^n$. If all functions $f$ and $h_i~(i=1,\ldots,m)$ are linear functions, problem \reff{prob1-1}--\reff{prob1-2} is a standard form linear programming problem (for examples, see \cite{NocWri99,wright97,ye}).
In this paper, we mainly focus on the nonlinear programs that at least one of functions $f$ and $h_i~(i=1,\ldots,m)$ is a nonlinear (and possibly nonconvex) function in problem \reff{prob1-1}--\reff{prob1-2}.
Our method can be easily extended to cope with nonlinear programs with general nonlinear inequality constraints (see section 6 for details).

There are already many efficient algorithms and several efficient solvers for nonlinear program \reff{prob1-1}--\reff{prob1-2}, among them is the state-of-the-art and well known solver LANCELOT (see \cite{ConGoT92}). Using the augmented Lagrangian function on equality constraints, Conn, Gould and Toint \cite{ConGoT92} solves the relaxed subproblem \bea \min_{x}\quad L_A(x,\lambda;\rho)\equiv f(x)-\lambda^Th(x)+\hf\rho\|h(x)\|^2\quad\hbox{s.t.}\quad x\ge 0, \label{cgtp}\eea
where $\la\in\Re^m$ is an estimate of the multiplier vector, $\rho>0$ is a penalty parameter. Both $\la$ and $\rho$ are held fixed during the solution of each subproblem and are updated adaptively in virtue of the convergence and feasibility of the approximate solution of the subproblem. Problem \reff{cgtp} is a nonlinear program with nonnegative constraints, and many algorithms in the literature can be used to solve this problem (see \cite{ConGoT88}).

Primal-dual interior-point methods have been demonstrated to be a class of very efficient methods for solving problem \reff{prob1-1}--\reff{prob1-2}. For example, for nonlinear programs, the readers can consult \cite{ByrGiN00,ByrHrN99,CheGol06,curtis12,CurGoR17,ForGil98,GerGil04,GoOrTo03,NocOzW12,ShaVan00,UlbUlV04,VanSha99,WacBie06} and the references there in. Generally, by requiring $x$ to be an interior-point, primal-dual interior-point methods solve the logarithmic-barrier subproblem \bea
\min_x\quad f(x)-\mu\sum_{j=1}^n\ln x_j\quad\hbox{s.t.}\quad h(x)=0 \label{lgbp}\eea
or its corresponding parametric Karush-Kuhn-Tucker (KKT) system, where $\mu>0$ is a barrier parameter which is held fixed when solving the subproblem \reff{lgbp} or its parametric KKT system. Different from subproblem \reff{cgtp} in the form, problem \reff{lgbp} is an equality constrained nonlinear program with logarithmic-barrier terms. Although all those effective algorithms for equality constrained nonlinear programming seem to be applicable to the subproblem, their convergence to a KKT point of the original problem may fail even for a well-posed problem (that is, a problem with a unique solution at which the second-order sufficient optimality conditions hold, see \cite{BenShV04,WacBie00}).

Improving the jamming difficulty (i.e., the failure of global convergence to a KKT point), the rapid convergence and the numerical performance of interior-point methods has been one of the main topics of the optimization research in recent years. For example, some warm-starting interior-point methods for linear programming have focused on relaxing the primal and dual interior-point limitations (see \cite{BS0,EngAnV09}) when the iterate is close to the solution. These methods were also extended to solve nonlinear programming in \cite{BS}. Numerical results in \cite{BS,EngAnV09} have shown that the warm-starting technique could improve the performance of interior-point methods for linear and nonlinear programming. Most recently, \cite{HHY19} investigated how the update of the barrier parameter affects the convergence of classic interior-point methods for convex and nonconvex optimization. Furthermore, \cite{HY18} proposed a one-phase interior-point method for nonconvex optimization with general inequality constraints, and showed that, by careful initialization and updates of the slack variables, the proposed method can be guaranteed to have more robust global convergence properties and will closely resemble successful algorithms from linear programming.

With the help of a logarithmic barrier augmented Lagrangian function, \cite{DLS17} proposed a bi-parametric primal-dual nonlinear system which corresponds to a KKT point and an infeasible stationary point of the original problem, respectively, as one of two parameters is zero. The method in \cite{DLS17} always generated interior-point iterates without any truncation of the step. Based on the equivalence of a positive relaxation problem to the logarithmic-barrier subproblem, \cite{LiuDai18} presented a globally convergent primal-dual interior-point relaxation method for nonlinear programs, which did not require any primal or dual iterate to be an interior point. The method has similarity to the warmstarting interior-point methods \cite{BS0,EngAnV09} and is different from most of the globally convergent interior-point methods in the literature.
Without assuming any regularity condition, the method either terminates at an approximate KKT point of
the original problem, an approximate infeasible stationary point, or an approximate singular stationary point of the original problem. The preliminary numerical results show that the algorithm is not only efficient for well-posed feasible problems, but also is applicable for some feasible problems without LICQ or MFCQ and some infeasible problems.

In this paper, we first prove that, under suitable conditions, any solution of a parametric equality constrained mini-max problem is a KKT point of the logarithmic-barrier subproblem. Based on this observation, we present a novel primal-dual interior-point relaxation method with iteratively updating barrier for nonlinear programs subject to general equality and nonnegative constraints. Our method is established on approximately solving a sequence of KKT systems of the parametric equality constrained mini-max subproblems, which avoids the requirement that any primal or dual iterate is an interior-point. The barrier parameter is updated with the iteration point as we did for linear programming, which is distinct from the newly proposed primal-dual interior-point relaxation method (see \cite{LiuDai18}) for nonlinear programming where the parameter is only updated in outer iterations when, for a fixed barrier, the inner iterations have found some approximate solutions of the logarithmic-barrier subproblems satisfying the given accuracy. In particular, our update for the barrier parameter is autonomous and iterative, allowing our method to
potentially avoid the possible difficulties caused by the inappropriate initial selection of the barrier parameter and to speed up convergence to the solution.

The method is easily extended for solving problems with general inequality constraints without incorporating any additional slack variables.
It has the potential to circumvent the jamming difficulty that appears with many interior-point methods for nonlinear programs and improve the ill conditioning of the existing primal-dual interior-point methods as the barrier parameter is small (see \cite{NocWri99}). Furthermore, a new smoothing approach, which is totally different from the techniques used in \cite{QSZ00}, is introduced to develop our relaxation method and promote convergence of the method. Under suitable conditions, it is proved that our method can be globally convergent and locally quadratically convergent to the KKT point of the original problem.
The preliminary numerical results on a well-posed problem for which many interior-point methods fail to find the minimizer and a set of test problems from the CUTEr collection show that our method is efficient.

Our paper is organized as follows. In section 2, we prove that the classic logarithmic-barrier subproblem can be equivalently converted into an equality constrained  mini-max problem. Based on this equivalence, we present the framework of our primal-dual interior-point relaxation method for nonlinear programs in section 3. In this section, we also figure out why our method can be expected to be efficient in improving the classic interior-point methods. We analyze and prove the global and local convergence results of our method for nonlinear programs in sections 4 and 5, respectively. Some preliminary numerical results on nonlinear programming test problems are reported in section 6. We conclude our paper in the last section.

Throughout the paper, we use standard notations from the literature. A letter with
subscript $k$ is related to the $k$th
iteration, the subscript $j$ indicates the $j$th component of a vector, and the subscript $kj$ is the $j$th
component of a vector at the $k$th iteration. All vectors are column vectors, and $z=(x,u)$ means $z=[x^T,\hspace{2pt}u^T]^T$. The expression
$\theta_k={O}(t_k)$ means that there exists a constant $M$
independent of $k$ such that $\rvert\theta_k\rvert\le M\rvert t_k\rvert$ for all $k$ large enough, and $\theta_k={o}(t_k)$ indicates that $\rvert \theta_k\rvert\le\epsilon_k\rvert t_k\rvert$ for all $k$ large enough with $\lim_{k\rightarrow	0}\epsilon_k=0$. If it is not specified, $I$ is an identity matrix whose order is either marked in the subscript or is clear in the context, and $\|\cdot\|$ is the Euclidean norm. Some unspecified notations may be
identified from the context.

\section{An equality constrained mini-max problem}\label{sec2}
Before presenting our main results, we review an equivalent problem of the logarithmic-barrier subproblem proposed in \cite{LiuDai18}.

For any given parameters $\mu\ge 0$ and $\rho>0$, and any $x\in\Re^n$ and $s\in\Re^{n}$, Liu and Dai \cite{LiuDai18} defined $z:\Re^{2n}\rightarrow \Re^n$, $z=z(x,s;\mu,\rho)$ and $y:\Re^{2n}\rightarrow \Re^n$, $y=y(x,s;\mu,\rho)$ by components to be functions on $(x,s)$ as follows, 
\bea
  z_j(x_j,s_j;\mu,\rho)\equiv\frac{1}{2\rho}\Big(\sqrt{(s_j-\rho x_j)^2+4\rho\mu}-(s_j-\rho x_j)\Big), \label{zydf1}\\
  y_j(x_j,s_j;\mu,\rho)\equiv\frac{1}{2\rho}\Big(\sqrt{(s_j-\rho x_j)^2+4\rho\mu}+(s_j-\rho x_j)\Big), \label{zydf}
\eea
where $j=1,\ldots,n$, $x\in\Re^n$ and $s\in\Re^{n}$ are variables\footnote{A little change is that both $z$ and $y$ are divided by $\rho$ in this paper.}. Based on definitions \reff{zydf1} and \reff{zydf},
Liu and Dai \cite{LiuDai18} proposed to solve an equivalent positive relaxation problem to the logarithmic-barrier subproblem \reff{lgbp} (see Theorem 2.3 of \cite{LiuDai18}) in the form
\begin{align}
\min_{x,s}  &~f(x)-\mu\sum_{j=1}^n \ln z_j(x,s;\mu,\rho) \label{ldprob1}\\
\st  &~h(x)=0, \label{ldprob2}\\
  &~z(x,s;\mu,\rho)-x=0.\label{ldprob3}
\end{align}

For convenience of readers and our subsequent discussions, we list some preliminary results in the following lemmas. These results have some similarities to Lemmas 2.1 and 2.2 and Theorem 2.3 of \cite{LiuDai18}.
\ble\label{lemzp} For given $\mu\ge 0$ and $\rho>0$, $z_j$ and $y_j$
are defined by \reff{zydf1} and \reff{zydf}. Then \\
(1) $z_j\ge 0$, $y_j\ge 0$, $z_j-x_j=y_j-(s_j/\rho)$, and $z_jy_j=\mu/\rho$; \\[5pt]
(2) $x_j\ge 0,\ s_j\ge 0,\ x_js_j=\mu$ if and only if $z_j-x_j=0$; \\[5pt]
(3) $z_j-x_j=\frac{\mu-x_js_j}{\rho(y_j+x_j)}$ and $\rho(z_j+y_j)=\sqrt{(s_j-\rho x_j)^2+4\rho\mu}$. \ele
\begin{proof}
Results (1) and (2) can be proved in the same way as Lemma 2.1 of Liu and Dai \cite{LiuDai18}. We are left to prove the result (3). Note that
\begin{align} z_j-x_j &= \frac{1}{2\rho}\Big(\sqrt{(s_j-\rho x_j)^2+4\rho\mu}-(s_j+\rho x_j)\Big) \nn\\
 &= \frac{2\mu-2x_js_j}{\sqrt{(s_j-\rho x_j)^2+4\rho\mu}+(s_j+\rho x_j)} \nn\\
 &= \frac{\mu-x_js_j}{\rho(y_j+x_j)}, \nn
 \end{align}
and the last equality in \refl{lemzp} (3) follows from the definitions \reff{zydf1} and \reff{zydf}.
All results are derived. 
\end{proof}

By \refl{lemzp}, we always have $\rho (y_j+x_j-z_j)=s_j$ and $\mu=\rho z_jy_j$. Moreover, it follows from \refl{lemzp} (3), $\rho(z-x)^T(y+z)=(n\mu-x^Ts)+\rho\|z-x\|^2$.

\ble\label{yp} Given $\mu>0$ and $\rho>0$. Let $z_j$ and $y_j$ be defined by \reff{zydf1} and \reff{zydf}. Then \\
(1) $z_j$ and $y_j$ are differentiable, respectively, on $x$ and $s$, and \bea
&&\na_xz_j=\frac{z_j}{z_j+y_j}e_j, \quad \na_xy_j=-\frac{y_j}{z_j+y_j}e_j, \label{20140327a}\\
&&\na_sz_j=-\frac{1}{\rho}\frac{z_j}{z_j+y_j}e_j, \quad
\na_sy_j=\frac{1}{\rho}\frac{y_j}{z_j+y_j}e_j, \label{20140327b} \eea
where $e_j\in\Re^n$ is the $j$-th coordinate vector; \\
(2) $z_j$ and $y_j$ are differentiable on $\mu$, and \bea
\frac{\partial z_j}{\partial \mu}=\frac{\partial y_j}{\partial \mu}=\frac{1}{\rho}\frac{1}{z_j+y_j}; \label{20190915a}\eea
(3) $z_j$ and $y_j$ are differentiable on $\rho$, and \bea
\frac{\partial z_j}{\partial \rho}=\frac{1}{\rho}\frac{z_j}{z_j+y_j}(x_j-z_j),\quad \frac{\partial y_j}{\partial \rho}=-\frac{1}{\rho}\frac{y_j}{z_j+y_j}(y_j+x_j). \nn\eea Thus, \bea
\frac{\partial (z_j-x_j)^2}{\partial \rho}=-\frac{2}{\rho}\frac{z_j}{z_j+y_j}(z_j-x_j)^2. \label{20190622a}\eea
\ele\begin{proof}  By the result (1) of Lemma 2.2 of Liu and Dai \cite{LiuDai18}, one has \bea
&&\na_x(\rho z_j)=\rho\frac{\rho z_j}{\rho z_j+\rho y_j}e_j, \quad \na_x(\rho y_j)=-\rho\frac{\rho y_j}{\rho z_j+\rho y_j}e_j, \nn\\
&&\na_s(\rho z_j)=-\frac{\rho z_j}{\rho z_j+\rho y_j}e_j, \quad
\na_s(\rho y_j)=\frac{\rho y_j}{\rho z_j+\rho y_j}e_j. \nn\eea
Thus, \reff{20140327a} and \reff{20140327b} follow immediately.

Due to  \bea
\frac{\partial(\rho z_j)}{\partial \mu}=\frac{\partial(\rho y_j)}{\partial \mu}=\hf\frac{4\rho}{2\sqrt{(s_j-\rho x_j)^2+4\rho\mu}}, \nn
\eea
the result \reff{20190915a} is derived from \refl{lemzp} (3).

Since $\rho(z_j+y_j)=\sqrt{(s_j-\rho x_j)^2+4\rho\mu}$ and $\rho(z_j-y_j)=\rho x_j-s_j$, one has
\bea  \frac{\partial \rho(z_j+y_j)}{\partial \rho}=\frac{1}{\rho}\frac{(\rho x_j-s_j)x_j+2\mu}{z_j+y_j}, \quad
\frac{\partial \rho(z_j-y_j)}{\partial \rho}=x_j. \nn\eea
Thus, \begin{align}
\frac{\partial z_j}{\partial \rho} &= \frac{1}{\rho}\left(\frac{\partial\rho z_j}{\partial \rho}-z_j\right) \nn\\
 &= \frac{1}{\rho}\left(\hf\Big(\frac{\partial \rho(z_j+y_j)}{\partial \rho}+\frac{\partial \rho(z_j-y_j)}{\partial \rho}\Big)-z_j\right) \nn\\
 &= \frac{1}{\rho}\left(\hf\Big(\frac{2\mu-x_j(s_j-\rho x_j)}{\rho(z_j+y_j)}+x_j\Big)-z_j\right) \nn\\
 &= \frac{1}{\rho}\frac{z_j}{z_j+y_j}(x_j-z_j), \nn\\
\frac{\partial y_j}{\partial \rho} &= \frac{\partial z_j}{\partial\rho}+\frac{1}{\rho}\left(z_j-y_j-x_j\right) \nn\\
 &= -\frac{1}{\rho}\frac{y_j}{z_j+y_j}(y_j+x_j), \nn
 \end{align}
\bea
\frac{\partial (z_j-x_j)^2}{\partial \rho}=2(z_j-x_j)\frac{\partial z_j}{\partial \rho}=-\frac{2}{\rho}\frac{z_j}{z_j+y_j}(z_j-x_j)^2. \nn\eea
This result implies that $\|z-x\|^2$ is a monotonically nonincreasing function on $\rho$. \end{proof} 

The following result is the foundation of development of the primal-dual interior-point relaxation method in \cite{LiuDai18}.
\ble Given $\mu>0$ and $\rho>0$. Let $(x^*,\la^*)$ be a KKT pair of the logarithmic-barrier subproblem
\reff{lgbp} and $(x^*,\la^*,s^*)$ satisfies its KKT system
\begin{align} 
&\na f(x^*)-\na h(x^*)\la^*-s^*=0, \label{bkkt11}\\
 & h(x^*)=0, \label{bkkt12}\\
 & x_j^*>0,\ s_j^*>0,\ x_j^*s_j^*=\mu,\ j=1,\ldots,n, \label{bkkt13}
\end{align}
where $\la^*\in\Re^{m}$ is the Lagrange multiplier vector. Then
$((x^*,s^*),(\la^*,s^*))$ is a KKT pair of the relaxation problem \reff{ldprob1}--\reff{ldprob3}.

Conversely, if $\mu>0$ and $\rho>0$, $((x^*,s^*),(\la^*,\nu^*))$ is a KKT pair of
problem \reff{ldprob1}--\reff{ldprob3}, where $\la^*\in\Re^{m}$ and
$\nu^*\in\Re^n$ are, respectively, the associated Lagrange multipliers of constraints
\reff{ldprob2} and \reff{ldprob3}, then $\nu^*=s^*$ and $(x^*,\la^*,s^*)$
satisfies the system \reff{bkkt11}--\reff{bkkt13}. Thus, $(x^*,\la^*)$ is a KKT
pair of the logarithmic-barrier subproblem \reff{lgbp}. \ele
\begin{proof} 
Please refer to the proof of Theorem 2.3 of \cite{LiuDai18}. 
\end{proof} 

Throughout the paper, we take $z$ and $y$ to be functions on $(x,s)$ dependent on parameters $(\mu,\rho)$. When it is thought to be clear in the context, we may ignore the variables and parameters in writing functions $z$ and $y$ for simplicity.

Now we consider the relaxation problem \reff{ldprob1}--\reff{ldprob3}. By incorporating the ``similar" augmented Lagrangian terms on constraints of \reff{ldprob3} into the objective function, and taking the maximum with respect to $s$, we obtain a particular mini-max problem
\bea\min_{x\in\{x\in\Re^n\rvert h(x)=0\}}\left\{f(x)+\sum_{j=1}^n\max_{s_j\in\Re} G(x_j,s_j;\mu,\rho)\right\}, \label{spp1} \eea
or its equivalent form \bea \min_{x\in\{x\in\Re^n\rvert h(x)=0\}}\max_{s\in\Re^n} F(x,s;\mu,\rho),  \nn\eea
where $F:\Re^{2n}\rightarrow\Re$, $F(x,s;\mu,\rho)\equiv f(x)+\sum_{j=1}^nG(x_j,s_j;\mu,\rho)$ and $G:\Re\rightarrow\Re$,
\begin{align*}
G(x_j,s_j;\mu,\rho)\equiv &-\mu\ln z_j(x_j,s_j;\mu,\rho)+s_j(z_j(x_j,s_j;\mu,\rho)-x_j)\\
&+\hf\rho\rvert z_j(x_j,s_j;\mu,\rho)-x_j\rvert^2.
\end{align*}
It should be noticed that the extra two terms $s^T(z(x,s;\mu,\rho)-x)+\hf\rho\|z(x,s;\mu,\rho)-x\|^2$ in $F(x,s;\mu,\rho)$ (comparing to \reff{ldprob1}) are not the usual augmented Lagrangian terms, since they definitely use the variables of $s$ of the function $z$ as the estimates of Lagrange multipliers, and take the parameter $\rho$ in $z$ as the penalty parameter. Moreover, the barrier parameter $\mu$ is used not only in the logarithmic-barrier terms but also in the other terms.

Using the previous preliminary results, we can derive some properties on $F(x,s;\mu,\rho)$.
\ble\label{bfprop} Given $\mu>0$ and $\rho>0$. Let $z=z(x,s;\mu,\rho)$ and $y=y(x,s;\mu,\rho)$ be defined by \reff{zydf1} and \reff{zydf}, $Z=\hbox{diag}(z)$, $Y=\hbox{diag}(y)$. \\
(1) If $f$ is twice differentiable, then $F$ is twice differentiable with respect to $x$ and $s$. Moreover,
\bea   
&\na_xF(x,s;\mu,\rho)=\na f(x)-\rho y,\quad &\na_x^2F(x,s;\mu,\rho)=\na^2 f(x)+\rho (Z+Y)^{-1}Y, \nn\\
  &\na_sF(x,s;\mu,\rho)=z-x,\quad &\na_s^2F(x,s;\mu,\rho)=-\frac{1}{\rho}(Z+Y)^{-1}Z. \nn 
  \eea
(2) Function $F(x,s;\mu,\rho)$ is a strictly concave function with respect to $s$, and $F(x,s;\mu,\rho)-f(x)$ is a strictly convex function with respect to $x$.\\
(3) There holds 
\bea
\frac{\partial F(x,s;\mu,\rho)}{\partial\rho}=\frac{(\rho-1)}{\rho}(z-x)^T(Z+Y)^{-1}Z(z-x). \nn\eea
\ele
\begin{proof}  Due to Lemmas \ref{lemzp} and \ref{yp}, one has the derivatives
\bea   
\frac{\partial G(x_j,s_j;\mu,\rho)}{\partial x_j}&=&\frac{-\mu-y_j(s_j+\rho z_j-\rho x_j)}{z_j+y_j}=-\rho y_j, \nn\\
  \frac{\partial G(x_j,s_j;\mu,\rho)}{\partial s_j}&=&z_j-x_j+\frac{\mu-z_js_j-\rho z_j(z_j-x_j)}{\rho(z_j+y_j)}=z_j-x_j. \nn
  \eea
Again by \refl{yp}, the second-order derivatives in (1) follow immediately.

The results in (2) are straightforward since $\na_s^2F(x,s;\mu,\rho)$ is always negative definite and $\na_x^2(F(x,s;\mu,\rho)-f(x))$ is always positive definite.

Note that $\mu=\rho z_jy_j$, \refl{yp} (3), \reff{20190622a}, and $
{\partial \ln z_j}/{\partial\rho}=z_j^{-1}{\partial z_j}/{\partial\rho},\
{\partial s_j(z_j-x_j)}/{\partial\rho}=s_j{\partial z_j}/{\partial\rho},$
the result (3) follows immediately due to $\rho(z_j-x_j)=\rho y_j-s_j$ and \bea
\frac{\partial G(x_j,s_j;\mu,\rho)}{\partial\rho}=\frac{\rho-1}{\rho}\frac{z_j}{z_j+y_j}(z_j-x_j)^2. \nn\eea
\end{proof} 

In the following, we prove our main result of this section, which is the foundation of our novel primal-dual interior-point relaxation method in this paper.
\begin{thm} \label{mrs} Let $\mu>0$ and $\rho>0$. The following two results can be obtained. \\
(1) The pair $(x^*,s^*)\in\Re^{n}\times\Re^n$ is a local solution of the mini-max problem \reff{spp1} if and only if $x^*>0$ is a local solution of the logarithmic-barrier subproblem \reff{lgbp} and $s^*_j=\mu/x^*_j$ for all $j=1,\cdots,n$. \\
(2) If $(x^*,s^*)\in\Re^{n}\times\Re^n$ is a local solution of the mini-max problem \reff{spp1} and $\na h(x^*)$ is of full column rank, then there exists a $\lambda^*\in\Re^m$ such that \bea
  \na f(x^*)-\na h(x^*)\la^*-s^*&=0, \label{mimkkt1-1}\\
  h(x^*)&=0, \label{mimkkt1-2}\\
  z^*-x^*&=0, \label{mimkkt1-3}
  \eea
where $z^*=z(x^*,s^*;\mu,\rho)$. Thus, $(x^*,\la^*)$ is a KKT pair of the logarithmic-barrier subproblem \reff{lgbp}. 
\end{thm} 
\begin{proof} 
(1) In light of \refl{bfprop}, for any $x_j>0$, $G(x_j,s_j;\mu,\rho)$ reaches its maximum at $s_j^*=\mu/x_j$ since $z_j(x_j,s_j^*;\mu,\rho)-x_j=0$.
If $x_j\le 0$, then $\frac{\partial G(x_j,s_j;\mu,\rho)}{\partial s_j}>0$, which means that $G(x_j,s_j;\mu,\rho)$ is strictly monotonically increasing to $\infty$ as $s_j\rightarrow \infty$.
Thus, \bea \max_{s_j\in\Re^n} G(x_j,s_j;\mu,\rho)=\left\{\ba{ll}
-\mu\ln x_j, & \hbox{if}\quad x_j>0; \\
\infty, & \hbox{otherwise,}\ea\right.\eea
and \bea \hbox{argmax}_{s_j\in\Re^n} G(x_j,s_j;\mu,\rho)=\left\{\ba{ll}
\mu/x_j, & \hbox{if}\quad x_j>0; \\
\infty, & \hbox{otherwise.}\ea\right.\eea
The result follows immediately from the above two equations.

(2) If $(x^*,s^*)$ is a solution of the mini-max problem \reff{spp1}, then $z^*-x^*=0$ by (1) and $x^*$ is a local solution of the subproblem \bea
\min_x  F(x,s^*;\mu,\rho) \label{20190623c}\\
\st  h(x)=0. \label{20190623d}\eea
Thus, if $\na h(x^*)$ is of full column rank, by the first-order necessary conditions of optimality (for example, see \cite{NocWri99,SunYua06}), there exists a $\lambda^*\in\Re^m$ such that $(x^*,\lambda^*)$ is a KKT pair of
subproblem \reff{20190623c}--\reff{20190623d}, i.e., there exists a $\lambda^*\in\Re^m$ such that \bea
  \na f(x^*)-\na h(x^*)\la^*-\rho y^*=0, \nn\\
  h(x^*)=0, \nn\\
  z^*-x^*=0, \nn\eea
where $y^*=y(x^*,s^*;\mu,\rho)$ and $z^*=z(x^*,s^*;\mu,\rho)$. Then the equations \reff{mimkkt1-1}--\reff{mimkkt1-3} are attained immediately since $z^*-x^*=0$ if and only if $y^*-s^*/\rho=0$ due to \refl{lemzp} (1).
\end{proof} 

Although the logarithmic-barrier subproblem \reff{lgbp}, its relaxation subproblem \reff{ldprob1}--\reff{ldprob3}, and the mini-max subproblem \reff{spp1} are equivalent in some sense, they provide us insightful views on the existing methods and possibilities for developing different and possibly robust methods for the original problem \reff{prob1-1}--\reff{prob1-2}. For example, by using the relaxation subproblem \reff{ldprob1}--\reff{ldprob3}, we can remove the interior-point restrictions on primal and dual variables in \cite{LiuDai18}. In this paper, we note that, $(x^*,s^*)$ is a solution of a mini-max subproblem if $x^*$ is a local solution of the logarithmic-barrier subproblem. Thus, the residual function on the system \reff{mimkkt1-1}--\reff{mimkkt1-3} is reasonable to be chosen as the merit function. In addition, by solving the system \reff{mimkkt1-1}--\reff{mimkkt1-3}, we are capable of improving the ill conditioning often observed during the final stages of the classic primal-dual algorithms based on solving the subproblem \reff{lgbp} or its corresponding KKT system (please refer to Section 3 for details).

As a special example, when $f$ and $h$ are linear functions, that is, program \reff{prob1-1}--\reff{prob1-2} is a linear programming problem, the mini-max problem is a particular saddle-point problem. The next result is a corollary of \reft{mrs}.
\bco Assume $\mu>0$ and $\rho>0$, $f$ and $h_i$ $(i=1,\ldots,m)$ are linear functions on $\Re^n$. The primal-dual pair $(x^*,s^*)$ is a solution of the mini-max problem \reff{spp1} if and only if there exists a $\lambda^*\in\Re^m$ such that $(x^*,\la^*)$ is a KKT pair of the logarithmic-barrier subproblem \reff{lgbp}. \eco

\section{A novel primal-dual interior-point relaxation method}

Based on solving the mini-max subproblem \reff{spp1}, we develop a novel primal-dual interior-point relaxation method for solving the nonlinear constrained optimization problem \reff{prob1-1}--\reff{prob1-2}. Since problem \reff{spp1} originates from the logarithmic-barrier subproblem, our method can be thought of as a variant of classic primal-dual interior-point methods. The method updates the barrier parameter $\mu$ in every iteration, which resembles some successful interior-point methods for linear and nonlinear programming (such as \cite{HY18,meh92,NocWaW09}), and is different from those based on the Fiacco-McCormick approach \cite{FM90} for nonlinear programming in which they often attempt to find an approximate solution for a fixed parameter $\mu$ in an inner algorithm and then reduce the barrier parameter $\mu$ by the residual of the solution in an outer algorithm. In particular, our update for the barrier parameter is autonomous and iterative, which makes our method capable of avoiding the possible difficulties caused by unappropriate initial selection of the barrier parameter and makes our method have the potential of speeding up the convergence to the solution.

Instead of solving the subproblem \reff{spp1} directly, we solve the associated system \reff{mimkkt1-1}--\reff{mimkkt1-3} and consider the extended system of equations of \reff{mimkkt1-1}--\reff{mimkkt1-3} in the form \bea
   \mu&=0, \label{mimkkt2-0}\\
   \na f(x)-\na h(x)\la-s&=0, \label{mimkkt2-1}\\
   h(x)&=0, \label{mimkkt2-2}\\
   z-x&=0, \label{mimkkt2-3}\eea
where $z=z(x,s;\mu,\rho)$ and $y=y(x,s;\mu,\rho)$ are functions on $x$ and $s$ defined by \reff{zydf1} and \reff{zydf}.
Distinct from our recent work \cite{DLS17, LiuDai18} and many interior-point methods for nonlinear programs, we also take $\mu$ as a variable in the system \reff{mimkkt2-0}--\reff{mimkkt2-3} instead of only a parameter in the system \reff{mimkkt1-1}--\reff{mimkkt1-3} so that $\mu$ is updated with the iteration point.
This approach has been used successfully in smoothing Newton methods for nonlinear complementarity problems and box constrained variational inequalities (see \cite{QSZ00}), where $\mu$ is a vector of smoothing parameters.
Note that, for $j=1,\ldots,n$, 
\bea 
z_j(x_j,s_j;0,\rho)&=&\frac{1}{2\rho}(\rvert s_j-\rho x_j\rvert-(s_j-\rho x_j))=\max\{0, x_j-s_j/\rho\}, \nn\\
y_j(x_j,s_j;0,\rho)&=&\frac{1}{2\rho}(\rvert s_j-\rho x_j\rvert+(s_j-\rho x_j))=\max\{0, s_j/\rho-x_j\}. \nn
\eea
Thus, for any $j=1,\ldots,n$, the equality $z_j=x_j$ implies that one has either $x_j=0$, $s_j\ge 0$, $\rho y_j=s_j$, or $x_j\ge 0$, $s_j=0$, $\rho y_j=0$. Therefore, any $(x^*, \la^*, s^*)\in\Re^n\times\Re^m\times\Re^n$ satisfying the extended system of equations \reff{mimkkt2-0}--\reff{mimkkt2-3} is a KKT triple of the original problem \reff{prob1-1}--\reff{prob1-2}.

Denote the residual function of the system \reff{mimkkt1-1}--\reff{mimkkt1-3} as follows, \bea \phi_{(\mu,\rho)}(x,\la,s)=\hf\|\na f(x)-\na h(x)\la-s\|^2+\hf\|h(x)\|^2+\hf\|z-x\|^2. \label{meritres}\eea
Using this notation, the system \reff{mimkkt2-0}--\reff{mimkkt2-3} can be further reformulated as
\bea
   \mu+\gamma\phi_{(\mu,\rho)}(x,\la,s)&=0, \label{mimkkt3-0}\\
   \na f(x)-\na h(x)\la-s&=0, \label{mimkkt3-1}\\
   h(x)&=0, \label{mimkkt3-2}\\
   z-x&=0, \label{mimkkt3-3}\eea
where $\mu$ is supposed to be nonnegative, $\phi_{(\mu,\rho)}(x,\la,s)$ is defined by \reff{meritres}, and $\gamma\in(0,1]$ is a given parameter.

In order to solve the system \reff{mimkkt2-0}--\reff{mimkkt2-3} efficiently, $\mu$ should not approach zero too quickly. Thus it is important to balance the reduction of $\mu$ and the associated KKT residual of the mini-max subproblem. The methods in \cite{QSZ00} were established on solving the system with elaborately constructed perturbation of the Newton system, and the residual function of the whole system was taken as the merit function. In contrast, instead of solving the system \reff{mimkkt2-0}--\reff{mimkkt2-3} directly, we develop our relaxation method by solving the reformulation \reff{mimkkt3-0}--\reff{mimkkt3-3} and promote convergence of our method by reducing the residual function $\phi_{(\mu,\rho)}(x,\la,s)$.

Suppose that $(x_k,\la_k,s_k)$ is the current primal and dual iterates, $\mu=\mu_k$ and $\rho=\rho_k$ are current values of the barrier and penalty parameters.
Let $r_k^d=\na f(x_k)-\na h(x_k)\la_k-s_k$, $r_k^e=z_k-x_k$, and $r_k^h=h(x_k)$ be the residuals of equations in \reff{mimkkt3-1}--\reff{mimkkt3-3} at iterate $k$.
Our proposed method generates the new value of parameter $\mu_{k+1}$ by \bea \mu_{k+1}=(1-\alp_k)\mu_k+\gamma\alp_k\phi_{(\mu_k,\rho_k)}(x_k,\la_k,s_k) \nn\eea
and the new primal and dual iterates by a line search procedure
$$x_{k+1}=x_k+\alpha_kd_{xk},\quad \la_{k+1}=\la_k+\alp_kd_{\la k},\quad s_{k+1}=x_k+\alpha_kd_{sk},$$
where $(d_{xk},d_{\la k},d_{sk})$ is the search direction derived from the Newton's equations of system \reff{mimkkt3-1}--\reff{mimkkt3-3}, and $\alpha_k\in (0,1]$ is the step-size. At iterate $(x_k,\la_k,s_k)$ with $\mu=\mu_k$ and $\rho=\rho_k$, $(d_{xk},d_{\la k},d_{sk})$ is derived from solving the linearized system with respect to $(x,\la,s)$ and $\mu$ as the following
\bea &&\left[\ba{ccc}
B_k & -\na h(x_k) & -I \\
\na h(x_k)^T & 0 & 0 \\
(Z_k+Y_k)^{-1}Y_k  & 0 & \frac{1}{\rho_k}(Z_k+Y_k)^{-1}Z_k\ea\right]\left[\ba{c}
d_x \\
d_{\la} \\
d_s \ea\right]\nn\\
&&=\left[\ba{c}
-r_k^d  \\
-r_k^h\\
r_k^e+\frac{1}{\rho_k}\Delta\mu_k(Z_k+Y_k)^{-1}e\ea\right], \label{dsysa}\eea
where the term on the variation $\Delta\mu$ of $\mu$ is moved to the right-hand-side of the linearized equation.
The preceding system can also be equivalently written as the linear system with a symmetric coefficient matrix in the form
\bea &&\left[\ba{ccc}
B_k+\rho_k(Z_k+Y_k)^{-1}Y_k & -\na h(x_k) & -(Z_k+Y_k)^{-1}Y_k \\
-\na h(x_k)^T & 0 & 0 \\
-(Z_k+Y_k)^{-1}Y_k & 0 & -\frac{1}{\rho_k}(Z_k+Y_k)^{-1}Z_k\ea\right]\left[\ba{c}
d_x \\
d_{\la}\\
d_s \ea\right]\nn\\
&&=-\left[\ba{c}
\hat r_k^d-\Delta\mu_k(Z_k+Y_k)^{-1}e  \\
-r_k^h\\
r_k^e+\frac{1}{\rho_k}\Delta\mu_k(Z_k+Y_k)^{-1}e\ea\right], \nn\eea
where $B_k$ is the Hessian of the Lagrangian $L(x,\la,s)=f(x)-\la^Th(x)-s^Tx$ or its approximation at $(x_k,\la_k,s_k)$, $z_k=z(x_k,s_k;\mu_k,\rho_k)$, $y_k=y(x_k,s_k;\mu_k,\rho_k)$, $Z_k=\diag(z_k)$, $Y_k=\diag(y_k)$, $\Delta\mu_k=-\mu_k+\gamma\phi_{(\mu_k,\rho_k)}(x_k,\la_k,s_k)$,
$\hat r_k^d=\na f(x_k)-\na h(x_k)\la_k-\rho_ky_k$.

Since we are facing a mini-max subproblem, taking the residual function $\phi_{(\mu,\rho)}(x,\la,s)$ defined by \reff{meritres}
as the merit function is a natural and reasonable selection.
The step-size $\alp_k$ is selected such that the value of $\phi_{(\mu,\rho)}(x,\la,s)$ is sufficiently decreased when the iterate moves from point $(x_k,\la_k,s_k)$ to $(x_{k+1},\la_{k+1},s_{k+1})$ and the barrier parameter varies from $\mu_k$ to $\mu_{k+1}$, while the penalty parameter $\rho_k$ holds fixed. Then $\rho_k$ is updated adaptively to $\rho_{k+1}$ such that $\rho_{k+1}\ge\rho_k$.

In the following, we describe our algorithm for problem \reff{prob1-1}--\reff{prob1-2}, in which the parameter $\mu$ is updated with the iteration point. In our algorithm, scalars $\gamma$, $\gamma_0$ and $\eta$ are parameters used to balance the reduction of $\mu_k$ and $\phi_{(\mu_k,\rho_k)}(x_k,\la_k,s_k)$. That is, for given $\gamma_0\in (0,1)$ and $\eta>1$, $\mu_k\in [\gamma_0\phi_{(\mu_k,\rho_k)}(x_k,\la_k,s_k), \eta\phi_{(\mu_k,\rho_k)}(x_k,\la_k,s_k)]$ is thought to be in a good balance and will be updated normally by the Newton's step; otherwise, it will be reduced provided it is larger or fixed if it is smaller before proceeding to a new iteration. The scalar $\gamma\in (0,\gamma_0]$ is a balance parameter introduced in (3.6). Scalars $\delta$ and $\tau$ are parameters necessary for Armijo's line search procedure in \reff{steprulea} and scalar $\sigma$ is a given factor for the update of the penalty parameter.
\begin{algorithm}
	\caption{A novel primal-dual interior-point relaxation method for problem \reff{prob1-1}--\reff{prob1-2}}
	\label{alg2}
	{\small \alglist
		\item[Given] $(x_0,\la_0,s_0)\in\Re^{n}\times\Re^m\times\Re^n$, $B_0\in\Re^{n\times n}$, $\mu_0>0$, $\rho_0>0$, $\eta>1$, $\gamma_0,\delta,\tau,\sigma\in(0,1)$. Evaluate $z_0$ and $y_0$ by \reff{zydf1} and \reff{zydf}, compute $\phi_{(\mu_0,\rho_0)}(x_0,\la_0,s_0)$. Given $\epsilon\in (0,\mu_0)$, set $k:=0$.
		
		Set $\ell:=0$, $\mu_{k,\ell}=\mu_k$.
		
		Step 0.1 While $\mu_{k,\ell}>\max\{\eta\phi_{(\mu_{k,\ell},\rho_k)}(x_k,\la_k,s_k),\epsilon\}$, set $\mu_{k,\ell+1}=\mu_{k,\ell}/\eta$;
		
		\hspace{1cm}evaluate $z_k$ and $y_k$ by \reff{zydf1} and \reff{zydf} with $\mu=\mu_{k,\ell+1}$, compute $\phi_{(\mu_{k,\ell+1},\rho_k)}(x_k,\la_k,s_k)$,
		
		\hspace{1cm}set $\ell=\ell+1$, end.
		
		Set $\mu_k=\mu_{k,\ell}$, $\gamma=\min\{\gamma_0, \mu_k/\phi_{(\mu_k,\rho_k)}(x_k,\la_k,s_k)\}$.

		\item[While] $\mu_k\le\epsilon$ and $\phi_{(\mu_k,\rho_k)}(x_k,\la_k,s_k)\le\epsilon$, stop the algorithm.
		
		\item[Step] 1. Calculate $\Delta\mu_k$ by $\Delta\mu_k=-\mu_k+\gamma\phi_{(\mu_k,\rho_k)}(x_k,\la_k,s_k)$.
		
		\item[Step] 2. Solve the linear system \reff{dsysa} to obtain $d_k\equiv(d_{xk},d_{\la k},d_{sk})$.
		
		\item[Step] 3. Select the step-size $\alp_k\in (0,1]$ to be the maximal number in $\{1,\delta,\delta^2,\ldots\}$ such that the inequality
		\bea
		\phi_{(\mu_{k}+\alp_k\Delta\mu_k,\rho_k)}&(x_{k}+\alp_kd_{xk},\la_{k}+\alp_kd_{\la k},s_{k}+\alp_kd_{sk})\nn\\
		&\le(1-2\tau\alpha_k)\phi_{(\mu_k,\rho_k)}(x_k,\la_k,s_k) \label{steprulea}
		\eea
		is satisfied.
		
		\item[Step] 4. Set $\mu_{k+1}=\mu_k+\alp_k\Delta\mu_k$, $x_{k+1}=x_k+\alpha_kd_{xk}$, $s_{k+1}=s_k+\alpha_kd_{sk}$, and $\la_{k+1}=\la_k+\alp_kd_{\la k}$.
		
		\item[Step] 5. Update $\rho_{k}$ to $\rho_{k+1}=\max\{\rho_k,\sigma\|s_{k+1}\|_{\infty}/\max(\|x_{k+1}\|,1)\}$.
		Evaluate by \reff{zydf1} and \reff{zydf} $$z_{k+1}=z(x_{k+1},s_{k+1};\mu_{k+1},\rho_{k+1})\ \hbox{and}\ y_{k+1}=y(x_{k+1},s_{k+1};\mu_{k+1},\rho_{k+1}),$$ compute $\phi_{(\mu_{k+1},\rho_{k+1})}(x_{k+1},\la_{k+1},s_{k+1})$.
		
		Set $\ell:=0$, $\mu_{k+1,\ell}=\mu_{k+1}$.
		
		Step 5.1 While $\mu_{k+1,\ell}>\max\{\eta\phi_{(\mu_{k+1,\ell},\rho_{k+1})}(x_{k+1},\la_{k+1},s_{k+1}),\epsilon\}$, set $\mu_{k+1,\ell+1}=\mu_{k+1,\ell}/\eta$;
		
		\hspace{1cm}evaluate $z_{k+1}$ and $y_{k+1}$ by \reff{zydf1} and \reff{zydf} with $\mu=\mu_{k+1,\ell+1}$,
		
		\hspace{1cm}compute $\phi_{(\mu_{k+1,\ell+1},\rho_{k+1})}(x_{k+1},\la_{k+1},s_{k+1})$, set $\ell=\ell+1$, end.
		
		Set $\mu_{k+1}=\mu_{k+1,\ell}$, $\gamma=\min\{\gamma_0, \mu_{k+1}/\phi_{(\mu_{k+1},\rho_{k+1})}(x_{k+1},\la_{k+1},s_{k+1})\}$.
		
		\item[Step] 6. Update $B_k$ to $B_{k+1}$, set $k:=k+1$.
		
		\item[End] (while)
		
		\eli}
	
\end{algorithm}

For \refal{alg2}, the initial point can be any point which is either an interior or other point. Our algorithm does not also require any primal or dual iterate to be interior during the iterative process, which is distinct from most of the classic interior-point methods.
Steps 0.1 and 5.1 are used to prevent $\mu_0$ and $\mu_{k+1}$ from being too large in comparison with the residuals of KKT system $\phi_{(\mu_{0},\rho_{0})}(x_{0},\la_{0},s_{0})$ and $\phi_{(\mu_{k+1},\rho_{k+1})}(x_{k+1},\la_{k+1},s_{k+1})$, respectively. If $\mu_{k+1}\le\max\{\eta\phi_{(\mu_{k+1},\rho_{k+1})}(x_{k+1},\la_{k+1},s_{k+1}),\epsilon\}$,
then one of the following three kinds of results will arise:

\quad(1) $\epsilon<\mu_{k+1}\le\eta\phi_{(\mu_{k+1},\rho_{k+1})}(x_{k+1},\la_{k+1},s_{k+1})$;

\quad(2) $\mu_{k+1}\le\epsilon\le\eta\phi_{(\mu_{k+1},\rho_{k+1})}(x_{k+1},\la_{k+1},s_{k+1})$;

\quad(3) $\mu_{k+1}\le\epsilon$ and $\phi_{(\mu_{k+1},\rho_{k+1})}(x_{k+1},\la_{k+1},s_{k+1})\le\epsilon/\eta<\epsilon$. \\
Note that, if the case (3) happens, \refal{alg2} will be terminated; otherwise, one will have either case (1) or case (2), and in both cases,
\bea\phi_{(\mu_{k+1},\rho_{k+1})}(x_{k+1},\la_{k+1},s_{k+1})\ge\epsilon/\eta. \label{20191230a}\eea
Moreover, for cases (1) and (2),
the parameter $\gamma$ is selected such that either $\gamma=\gamma_0$ and $\mu_{k+1}>\gamma_0\phi_{(\mu_{k+1},\rho_{k+1})}(x_{k+1},\la_{k+1},s_{k+1})$ or $\Delta\mu_{k+1}=0$. If $\mu_{k+1}>\gamma_0\phi_{(\mu_{k+1},\rho_{k+1})}(x_{k+1},\la_{k+1},s_{k+1})$, then \bea
\mu_{k+2,0}=(1-\alp_{k+1})\mu_{k+1}+\alp_{k+1}\gamma_0\phi_{(\mu_{k+1},\rho_{k+1})}(x_{k+1},\la_{k+1},s_{k+1})<\mu_{k+1}\quad \label{20191230b}\eea
and \bea\mu_{k+2,0}>\gamma_0\phi_{(\mu_{k+1},\rho_{k+1})}(x_{k+1},\la_{k+1},s_{k+1})\ge(\gamma_0/\eta)\epsilon; \label{20191230c}\eea
otherwise, $$\mu_{k+1}=\gamma\phi_{(\mu_{k+1},\rho_{k+1})}(x_{k+1},\la_{k+1},s_{k+1})\le\gamma_0\phi_{(\mu_{k+1},\rho_{k+1})}(x_{k+1},\la_{k+1},s_{k+1}),$$ 
$\mu_{k+1}$ is viewed as to be too small in comparison with $\phi_{(\mu_{k+1},\rho_{k+1})}(x_{k+1},\la_{k+1},s_{k+1})$ and set $\mu_{k+2,0}=\mu_{k+1}$. Thus, there is always $(\gamma_0/\eta)\epsilon\le\mu_{k+1,0}\le\mu_{k}$ for all $k\ge 0$.

In order to have a deep understanding on the significance of \refal{alg2}, let us consider its application to the linear programs with the standard form
\bea \min\ c^Tx \quad\st\ Ax=b,\ x\ge 0. \eea
Corresponding to the original problem \reff{prob1-1}--\reff{prob1-2}, $f(x)=c^Tx$, $h(x)=Ax-b$. In this case, $\na f(x)=c$ and $\na h(x)=A^T$. Without loss of generality, we suppose that $A$ has full row rank. Since the Lagrangian Hessian is null, \reff{dsysa} is reduced to the following system \bea &&\left[\ba{ccc}
0  & -A^T & -I \\
A & 0 & 0 \\
(Z_k+Y_k)^{-1}Y_k & 0 & \frac{1}{\rho_k}(Z_k+Y_k)^{-1}Z_k \ea\right]\left[\ba{c}
d_x \\
d_{\la} \\
d_s\ea\right]\nn\\
&&=\left[\ba{l}
-(c-A^T\la_k-s_k) \\
-(Ax_k-b) \\
(z_k-x_k)+\frac{1}{\rho_k}\Delta\mu_k(Z_k+Y_k)^{-1}e \ea\right], \nn\eea
which, due to \refl{lemzp} (3), can be further written as 
\bea \left[\ba{ccc}
0  & A^T & I \\
A & 0 & 0 \\
\rho_k Y_k & 0 & Z_k \ea\right]\left[\ba{c}
d_x \\
d_{\la} \\
d_s\ea\right]&=&\left[\ba{l}
c-A^T\la_k-s_k \\
b-Ax_k \\
\mu_k e-X_kS_ke \ea\right]+\nn\\
&&\left[\ba{c}
0 \\
0 \\
\rho_k(Z_k-X_k)(z_k-x_k)+\Delta\mu_ke \ea\right], \quad\label{dsysa1}\eea
where the minus signs in the first row are changed by left multiplying a negative identity matrix and the last row in the system is derived by left multiplying $\rho_k(Z_k+Y_k)$, respectively, on both sides of the equations.

Comparing with the system in classic primal-dual interior-point methods for linear programming (for example, see (14.12) of Nocedal and Wright \cite{NocWri99}), our system \reff{dsysa1} is different in that both $S_k$ and $X_k$ in the last row of the Jacobian have been substituted with $\rho Y_k$ and $Z_k$ and the associated right-hand-side term has also been changed (i.e., some additional correction terms have been incorporated). As we will note from what follows, these changes make our method capable of improving the ill conditioning of primal-dual interior-point methods for linear programming.

Note that \reff{dsysa1} can be formalized as
\begin{eqnarray}
&&(AY_k^{-1}Z_kA^T)d_{\lambda}=\rho_k (b-Ax_k)+AY_k^{-1}Z_k(c-A^T\la_k-s_k) \nn\\
&&\hspace{3cm}-\rho_k A(I+Y_k^{-1}Z_k)(z_k-x_k)-\Delta\mu_k AY_k^{-1}e, \nonumber\\[2pt]
&&d_s=(c-A^T\la_k-s_k)-A^Td_{\lambda}, \nonumber\\[2pt]
&&d_x=(I+Y_k^{-1}Z_k)(z_k-x_k)+\frac{1}{\rho_k}(\Delta\mu_kY_k^{-1}e-Y_k^{-1}Z_kd_s). \nonumber
\end{eqnarray}
Due to $\rho_k Y_kZ_k=\mu_k I$, one has $Y_k^{-1}=(\rho_k/\mu_k)Z_k$, and
\begin{eqnarray}
&&(AZ_k^2A^T)d_{\lambda}=\mu_k (b-Ax_k)+AZ_k^2(c-A^T\la_k-s_k) \nn\\
&&\hspace{2.5cm}-A(\mu_kI+\rho_k Z_k^2)(z_k-x_k)-\Delta\mu_k AZ_ke, \label{linsys1-1}\\[2pt]
&&d_s=(c-A^T\la_k-s_k)-A^Td_{\lambda}, \label{linsys1-2}\\[2pt]
&&d_x=(I+\frac{\rho_k}{\mu_k} Z_k^2)(z_k-x_k)+\frac{1}{\mu_k}(\Delta\mu_kZ_ke-Z_k^2d_s). \label{linsys1-3}
\end{eqnarray}
If $z_k\rightarrow  x^*$ and $\mu_k\rightarrow  0$ as $k\rightarrow \infty$, where $x^*$ is an optimal solution of the nondegenerate linear program, then $AZ_k$ should be of full rank and \reff{linsys1-1} is capable of escaping from the ill conditioning trap often observed during the final stages of the existing primal-dual algorithms for linear programming (see, for example, page 409 of \cite{NocWri99}). One may note that \reff{linsys1-3} could be possibly numerically difficult as $\mu_k\rightarrow  0$. However, in contrast to the implicit trap of the existing primal-dual algorithms, this difficulty of \reff{linsys1-3} is explicit and singlet. Theoretically, under suitable conditions, we can prove that, for all $j=1,\ldots,n$, $\frac{1}{\mu_k}z_{kj}$ is bounded away from zero (see \refl{boundness} for details).

Subsequently, we will show that \refal{alg2} is well-defined. Firstly, it is easy to note that Steps 0.1 and 5.1 will always be terminated finitely for any given $\epsilon>0$.
\ble There always holds $\mu_{k}\ge(\gamma_0/\eta)\epsilon$ for all $k\ge 0$. \ele
\begin{proof} 
We firstly prove that, if $\mu_{k+1,0}>\max\{\eta\phi_{(\mu_{k+1,0},\rho_{k+1})}(x_{k+1},\la_{k+1},s_{k+1}),\epsilon\}$, then \bea\mu_{k+1}\ge\epsilon/\eta. \label{20200405a}\eea
By Step 5.1, $\mu_{k+1}=\mu_{k+1,\ell}$ for some $\ell\ge 1$. Thus, $\mu_{k+1,\ell-1}>\epsilon$ and $\mu_{k+1,\ell}=\mu_{k+1,\ell-1}/\eta$, which implies $\mu_{k+1,\ell}>\epsilon/\eta$. If $\mu_{k+1,0}\le\max\{\eta\phi_{(\mu_{k+1,0},\rho_{k+1})}(x_{k+1},\la_{k+1},s_{k+1}),\epsilon\}$, $\mu_{k+1}=\mu_{k+1,0}$.

We have already known that $\mu_{k+1,0}\ge(\gamma_0/\eta)\epsilon$. Note that $\mu_0>\epsilon$, the result follows immediately from \reff{20200405a}.
\end{proof} 

In view of \reff{zydf1} and \reff{zydf}, $\mu_k>0$ implies $y_k>0$ and $z_k>0$. The following result asserts that the linear system \reff{dsysa} has a unique solution. \ble\label{lem3.2} Let $x_k$ be the current iterate generated by \refal{alg2}.
If $\na h(x_k)$ has full column rank and $v^T(B_k+\rho_kZ_k^{-1}Y_k)v>0$ for all $v\ne 0\in\Re^n$ with $\na h(x_k)^Tv=0$, then the coefficient matrix of the linear system \reff{dsysa} is nonsingular.
\ele\begin{proof}  In order to obtain our desired result, we need prove that the system of equations \bea
& B_kd_x-\na h(x_k)d_{\la}-d_s=0, & \label{20191111a}\\
& \na h(x_k)^Td_x=0, & \label{20191111b}\\
& \rho_k Y_kd_x+Z_kd_s=0 \label{20191111c}\eea
has only zero solution. Left-multiplying $d_x^T$ on the two-sides of \reff{20191111a}, one has $d_x^TB_kd_x=d_x^Td_s$ due to \reff{20191111b}. Thus, by \reff{20191111c}, \bea d_x^T(B_k+\rho_kZ_k^{-1}Y_k)d_x=0. \label{20191111d} \eea
Note that the conditions of the lemma suggest $d_x^T(B_k+\rho_kZ_k^{-1}Y_k)d_x>0$ for all $d_x\ne 0$ satisfying \reff{20191111b}, thus $d_x=0$. Therefore, $d_s=0$ and $\na h(x_k)d_{\la}=0$ due to the last and the first equations of the preceding system. Since $\na h(x_k)$ has full column rank, the equation $\na h(x_k)d_{\la}=0$ implies $d_{\la}=0$. Hence, our proof is completed.
\end{proof} 

If \refal{alg2} does not terminate at $x_k$, then $\phi_{(\mu_k,\rho_k)}(x_k,\la_k,s_k)\ge\epsilon/\eta>0$ due to \reff{20191230a}.
This fact shows that there will be $(\De\mu_k,d_k)\ne 0$ for all $k\ge 0$. Otherwise, by \refl{lem3.2}, the right-hand-side of \reff{dsysa} will be zero for some integer $k$, which implies $r_k^e=r_k^d=r_k^h=0$. Thus, $\phi_{(\mu_k,\rho_k)}(x_k,\la_k,s_k)=0$, a contradiction to \reff{20191230a}. The next result shows that, at the $k$-th iteration, a new iterate can be generated, thus \refal{alg2} is well-defined.
\ble Suppose that $f: \Re^n\rightarrow \Re$ and $h: \Re^n\rightarrow \Re^m$ are twice continuously differentiable on $\Re^n$. There always exists an $\alp_k\in(0,1]$ such that \reff{steprulea} holds.
\ele
\begin{proof}  The supposition implies that $\phi_{(\mu,\rho_k)}(x,\la,s)$ is differentiable with respect to $(\mu,x,\la,s)$, thus it is directionally differentiable. Due to \reff{dsysa}, its directional derivative along $(\Delta\mu_k,d_{k})$ at $(x_k,\la_k,s_k)$ with $\mu=\mu_k$ is 
	\bea
  &&\phi'_{(\mu_k,\rho_k)}(x_k,\la_k,s_k;\Delta\mu_k,d_{k}) \nn\\
  &=&\left(\frac{D\phi_{(\mu,\rho_k)}(x,\la,s)}{D\mu}\quad \nabla_{(x,\la,s)}\phi_{(\mu,\rho_k)}(x,\la,s)^T\right)\rvert_{(\mu,x,\la,s)=(\mu_k,x_k,\la_k,s_k)}\left(\ba{c}
\Delta\mu_k\\
d_k\ea\right)\nn\\
  &=&-2\phi_{(\mu_k,\rho_k)}(x_k,\la_k,s_k). \eea

The Taylor's expansion of $\phi_{(\mu_{k}+\alp\Delta\mu_k,\rho)}(x_{k}+\alp d_{xk},s_{k}+\alp d_{\la k},\la_{k}+\alp d_{sk})$ regarding $\alp$ at $\alp=0$ shows that \bea   &&\phi_{(\mu_{k}+\alp\Delta\mu_k,\rho)}(x_{k}+\alp d_{xk},s_{k}+\alp d_{\la k},\la_{k}+\alp d_{sk})-\phi_{(\mu_k,\rho_k)}(x_k,\la_k,s_k) \nn\\
  &=&\alp\phi'_{(\mu_k,\rho_k)}(x_k,\la_k,s_k;\Delta\mu_k,d_{k})+o(\alp) \nn\\
  &=&-2\tau\alpha\phi_{(\mu_k,\rho_k)}(x_k,\la_k,s_k)-2(1-\tau)\alp\phi_{(\mu_k,\rho_k)}(x_k,\la_k,s_k)+o(\alp). \label{meritp}\eea
Thus, \reff{steprulea} holds for all sufficiently small $\alp>0$ since $\tau<1$ and $\phi_{(\mu_k,\rho_k)}(x_k,\la_k,s_k)>0$.
\end{proof} 

The preceding result suggests that sequences $\{(x_k,\la_k,s_k)\}$ and $\{\mu_k\}$, $\{\rho_k\}$ will be derived from \refal{alg2} before the terminating condition is satisfied.
Moreover, \reff{20191230b} has shown that the barrier sequence $\{\mu_k\}$ is monotonically nonincreasing.
It will be proved that the sequence of merit function values $\{\phi_{(\mu_k,\rho_k)}(x_k,\la_k,s_k)\}$ is monotonically decreasing.
\ble\label{lemrho0} Let $z_{k+1}(\rho)=z(x_{k+1},s_{k+1};\mu_{k+1},\rho)$ and $\hat z_{k+1}=z_{k+1}(\rho_k)$. Suppose that $\|\hat z_{k+1}-x_{k+1}\|\ne 0$ and $\phi_{(\mu_{k+1},\rho_k)}(x_{k+1},\la_{k+1},s_{k+1})\le\phi_{(\mu_k,\rho_k)}(x_k,\la_k,s_k).$ If $\rho_{k+1}\ge\rho_k>0$, one has
$$\phi_{(\mu_{k+1},\rho_{k+1})}(x_{k+1},\la_{k+1},s_{k+1})\le\phi_{(\mu_k,\rho_k)}(x_k,\la_k,s_k).$$
\ele\begin{proof}  Note that \bea   &&\frac{D\phi_{(\mu_{k+1},\rho)}(x_{k+1},\la_{k+1},s_{k+1})}{D\rho}\rvert_{\rho=\rho_k} \nn\\
  &&=-\frac{1}{\rho_k}(\hat z_{k+1}-x_{k+1})^T(\hat Z_{k+1}+\hat Y_{k+1})^{-1}\hat Z_{k+1}(\hat z_{k+1}-x_{k+1}) \nn\\
  &&<0, \nn\eea where $\hat Z_{k+1}=\diag(\hat z_{k+1})$, $\hat Y_{k+1}=\diag(\hat y_{k+1})$ with $\hat y_{k+1}=z(x_{k+1},s_{k+1};\mu_{k+1},\rho_k)$. The above equation shows that $\phi_{(\mu_{k+1},\rho)}(x_{k+1},\la_{k+1},s_{k+1})$ is a monotonically decreasing function on $\rho$ over $\rho>0$, which implies the desired result.
\end{proof} 

By \refal{alg2}, the sequence $\{\rho_k\}$ of penalty parameters is a monotonically nondecreasing sequence. The following result follows from Steps 0.1 and 5.1 immediately.
\ble\label{lemmu0} There hold $$0<\mu_{k+1}\le\mu_{k}\le\mu_0\quad\hbox{and}\quad
\mu_{k}\le\eta\phi_{(\mu_k,\rho_k)}(x_k,\la_k,s_k)$$ for all $k>0$.\ele\begin{proof} 
The result follows from \reff{20191230b} and $\mu_{k+2}\le\mu_{k+2,0}$ immediately.\end{proof}

\section{Global convergence}

For global and local convergence analysis, we set $\epsilon=0$. In this situation, \refal{alg2} may have infinite loop in either Step 0.1 for the initial iteration $k=0$ or in Step 5.1 for some iteration $k>0$. In any of these two trivial cases, one will have $\ell\rightarrow \infty$, $\lim_{\ell\rightarrow \infty}\mu_{k,\ell}=0$ and $\lim_{\ell\rightarrow \infty}\phi_{(\mu_{k,\ell},\rho_{\ell})}(x_{k},\la_{k},s_{k})=0$, thus $(x_k,\la_k,s_k)$ is a KKT triple of the problem \reff{prob1-1}--\reff{prob1-2}. Otherwise, \refal{alg2} will generate an infinite sequence of vectors $\{(x_k,\la_k,s_k)\}$. We consider this nontrivial case and prove in this section that, under suitable assumptions, there are some cluster points of the iterative sequence $\{(x_k,\la_k,s_k)\}$ which will be KKT triples of the problem \reff{prob1-1}--\reff{prob1-2}, i.e., the cluster points together with $\mu^*=0$ are solutions of the system of equations \reff{mimkkt2-0}--\reff{mimkkt2-3}.

We need the following blanket assumptions for our global convergence analysis.
\bas\label{ass1}\ \\
(1) The functions $f$ and $h_i~(i\in{\cal I})$ are twice continuously
differentiable on $\Re^n$; \\
(2) The iterative sequence $\{x_k\}$ is in an open bounded set of $\Re^n$;\\
(3) The sequence $\{B_k\}$ is bounded, and for all $k\ge 0$ and all $d_x\in\Re^n: d_x\ne 0, \na h(x_k)^Td_x=0$, $d_x^T(B_k+\rho_kZ_k^{-1}Y_k)d_x\ge\chi\|d_x\|^2,$ where $\chi>0$ is a constant;\\
(4) For all $k\ge 0$, $\na h(x_k)$ has full column rank.
\eas

The above assumptions are commonly used in global convergence analysis for nonlinear programs. Some milder assumptions can be used by incorporating some additional optimization techniques, such as the null-space technology (see \cite{BurCuW14,byrd,ByrGiN00,LiuSun01,LiuYua07}) for weakening \refa{ass1} (3) and (4), and the line search procedure without using a penalty function or a filter (see \cite{GouToi07,LiuYua08}) for replacing \refa{ass1} (2) on the requirement of the boundedness of the iterative sequence by some assumptions on bounded level sets. For simplicity of statement, we leave these concerns outside our scope. The following lemma shows that some related sequences are bounded.

\ble\label{boundness}
Under \refa{ass1}, $\{z_k\}$ is bounded and $\{s_k\}$ is bounded below. Furthermore, if
$[\na h(x_k)\quad I_{{\cal A}_k}]$ has full column rank for all $k$, where ${\cal A}_k=\{j\in\{1,2,\ldots,n\}\rvert{s_{kj}}\ne 0\}$, $I_{{\cal A}_k}\in\Re^{n\times \rvert{\cal A}_k\rvert}$ is a submatrix of $I_n$ with indices of the columns in ${\cal A}_k$, then
$\rho_k$ keeps constant after a finite number of iterations, $\{y_k\}$, $\{s_k\}$ and $\{\la_k\}$ are bounded, and there exists a scalar $\hat\tau>0$ such that, for $j=1,\ldots,n,$
\bea y_{kj}\ge\hat\tau\mu_k,\ z_{kj}\ge\hat\tau\mu_k.\nn\eea
\ele\begin{proof} 
Note that $z_k\ge 0$ for all $k\ge 0$ and
$$\phi_{(\mu_{k+1},\rho_{k+1})}(x_{k+1},\la_{k+1},s_{k+1})\le\phi_{(\mu_k,\rho_{k})}(x_k,\la_k,s_k)\le\ldots\le\phi_{(\mu_0,\rho_0)}(x_0,\la_0,s_0).$$
By the definition \reff{meritres} of $\phi_{(\mu,\rho)}(x,\la,s)$, one has
\bea \hf\|z_k-x_k\|^2\le\phi_{(\mu_0,\rho_0)}(x_0,\la_0,s_0), \nn\eea
which together with \refa{ass1} (2) implies that $\{z_k\}$ is bounded. Thus, due to \reff{zydf1}, for every $j=1,\ldots,n$,
$\sqrt{(s_{kj}/\rho_k-x_{kj})^2+4\mu_k/\rho_k}-(s_{kj}/\rho_k-x_{kj})$ is bounded. That is, $s_{kj}/\rho_k\not\rightarrow-\infty$ as $k\rightarrow \infty$, which implies that $\{s_k\}$ is bounded below.

Note that \bea\hf\|\na f(x_k)-\na h(x_k)\la_k-s_k\|^2\le\phi_{(\mu_0,\rho_0)}(x_0,s_0,\la_0).\label{rev090301}\eea
If there is a subsequence $\{s_{k_i}\}$ such that $\|s_{k_i}\|_{\infty}\rightarrow \infty$ as $k_i\rightarrow \infty$, then, due to \reff{rev090301}, one should have $\|\la_{k_i}\|_{\infty}\rightarrow \infty$ as $k_i\rightarrow \infty$. Divide by $\|(\la_k, s_k)\|_{\infty}$ and take the limit on the two sides of \reff{rev090301} as $k_i\rightarrow \infty$, it follows \bea \lim_{k_i\rightarrow \infty}\left\|\na h(x_{k_i})\frac{\la_{k_i}}{\|(\la_{k_i}, s_{k_i})\|_{\infty}}+\frac{s_{k_i}}{\|s_{k_i}\|}\right\|=0, \eea
which contradicts the condition that $[\na h(x_k)\quad I_{{\cal A}_k}]$ is of full column rank. The contradiction shows that $\{s_k\}$ and $\{\la_k\}$  are bounded. Furthermore, the update rule of $\rho_k$ implies that $\{\rho_k\}$ is bounded above. Thus, by \reff{zydf}, $\{y_k\}$ is bounded.

The relation $\rho_ky_{kj}z_{kj}=\mu_k$ together with that facts that both $\{y_k\}$ and $\{z_k\}$ are bounded implies the desired inequalities. \end{proof} 

The preceding results show that, under suitable conditions, $\rho_k$ will keep constant after a finite number of iterations. In other words, there exists a scalar $\rho^*>0$, such that $\rho_k=\rho^*$ for all sufficiently large $k$. In this situation, the sequence $\{\phi_{(\mu_k,\rho_k)}(x_k,\la_k,s_k)\}$ and the second derivatives of $\phi_{(\mu,\rho)}(x,\la,s)$ for all iterates are bounded.
In the following, we prove that there holds $\mu_k\rightarrow  0$ and $\phi_{(\mu_k,\rho_k)}(x_k,\la_k,s_k)\rightarrow  0$.
\ble\label{lemstepsize}
Under \refa{ass1}, suppose that $\rho_k=\rho^*$ for all sufficiently large $k$, where $\rho^*>0$ is a scalar. If $\mu_k\le\eta\phi_{(\mu_k,\rho_k)}(x_k,\la_k,s_k)$ for all sufficiently large $k$, then \bea
\lim_{k\rightarrow \infty}\phi_{(\mu_k,\rho_k)}(x_k,\la_k,s_k)=0\quad {\hbox{and}}\quad \lim_{k\rightarrow \infty}\mu_k=0. \nn\eea
\ele\begin{proof} 
Note that $\{\phi_{(\mu_k,\rho_k)}(x_k,\la_k,s_k)\}$ is a monotonically nonincreasing sequence. Thus, by the boundedness of $\{\phi_{(\mu_k,\rho_k)}(x_k,\la_k,s_k)\}$, there is a scalar $\phi^*\ge 0$ such that $$\lim_{k\rightarrow \infty}\phi_{(\mu_k,\rho_k)}(x_k,\la_k,s_k)=\phi^*,\quad \lim_{k\rightarrow \infty}\alp_k\phi_{(\mu_k,\rho_k)}(x_k,\la_k,s_k)=0.$$

We prove the result by contradiction. Assume that $\phi^*>0$. Then the preceding equations imply $\lim_{k\rightarrow \infty}\alp_k=0$ and $\liminf_{k\rightarrow \infty}\mu_k>0$ since $\mu_k$ keeps constant provided $\mu_k\le\gamma_0\phi_{(\mu_k,\rho_k)}(x_k,\la_k,s_k)$.
Hence, by \refl{boundness}, $z_k$ and $y_k$ are bounded away from zero. Similar to \refl{lem3.2}, we can prove that the matrix \bea \left(\ba{ccc}
B_k  & -A_k^T & -I \\
A_k & 0 & 0 \\
\rho_k Y_k & 0 & Z_k\ea\right) \nn\eea
is nonsingular for all $k$, where $A_k=\na h(x_k)^T$. Therefore, $\|d_k\|$ is bounded. In this case
\refa{ass1} asserts that $\alp_k$ is bounded away from zero since, by \reff{meritp}, \bea
  &\phi_{(\mu_{k}+\alp\Delta\mu_k,\rho_k)}&(x_{k}+\alp d_{xk},\la_{k}+\alp d_{\la k},s_{k}+\alp d_{sk})-(1-2\tau\alp)\phi_{(\mu_k,\rho_k)}(x_k,\la_k,s_k) \nn\\
  &&=-2(1-\tau)\alp\phi_{(\mu_k,\rho_k)}(x_k,\la_k,s_k)+o(\alp) \nn\\
  &&\le-2(1-\tau)\phi^*\alp+o(\alp), \nn\eea
which suggests that there exists an $\alp^*\in(0,1)$ such that \reff{steprulea} holds for all $\alp\in(0,\alp^*]$. It is contrary to $\lim_{k\rightarrow \infty}\alp_k=0$. This contradiction shows $\phi^*=0$. The desired results are obtained accordingly.
\end{proof} 

Now we are ready for presenting our global convergence results on \refal{alg2}.
\begin{thm} \label{mainglob}
Under \refa{ass1}, suppose that $\rho_k=\rho^*$ for all sufficiently large $k$, where $\rho^*>0$ is a scalar. Then one of the following three cases will arise. \\
(1) For all sufficiently large $k$, $\mu_k\le\eta\phi_{(\mu_k,\rho_k)}(x_k,\la_k,s_k)$. In this case, $\phi_{(\mu_k,\rho_k)}(x_k,\la_k,s_k)\rightarrow  0$ and $\mu_k\rightarrow  0$ as $k\rightarrow \infty$. That is, every cluster point of sequence $\{(x_k,\la_k,s_k)\}$ is a KKT triple of the original problem. \\
(2) For some iteration $k\ge 0$, $\mu_k>\eta\phi_{(\mu_k,\rho_k)}(x_k,\la_k,s_k)$, either Step 0.1 or Step 5.1 of \refal{alg2} has an infinite loop, $\lim_{\ell\rightarrow  0}\mu_{k,\ell}=0$ and $\lim_{\ell\rightarrow  0}\phi_{(\mu_{k,\ell},\rho_k)}(x_k,\la_k,s_k)=0$, i.e., $(x_k,\la_k,s_k)$ is a KKT triple of the original problem. \\
(3) Both Step 0.1 and Step 5.1 of \refal{alg2} have finite loops and Step 5.1 of \refal{alg2} is started over infinitely many times. Then
$ \lim_{k\rightarrow \infty}\mu_{k}=0$, and
there is an infinite subsequence $\{(x_{k_i},\la_{k_i},s_{k_i})\}$ of sequence $\{(x_k,\la_k,s_k)\}$ such that
\bea \lim_{i\rightarrow \infty}\phi_{(\mu_{k_i},\rho_{k_i})}(x_{k_i},\la_{k_i},s_{k_i})=0. \nn\eea
That is, there is a cluster point of sequence $\{(x_k,\la_k,s_k)\}$ is a KKT triple of the original problem.
\end{thm} \begin{proof}  The result in case (1) has been obtained in the preceding \refl{lemstepsize}. In case (2), let $\mu_{k,0}=\mu_k$ and $\mu_{k,\ell}=\mu_{k,\ell-1}/\eta$, where $\ell=1,2,\ldots$ is the number of the cycle of while in Step 5.1 of \refal{alg2}. Thus, $\lim_{\ell\rightarrow \infty}\mu_{k,\ell}=0$ and $\lim_{\ell\rightarrow \infty}\mu_{k,\ell}\ge\lim_{\ell\rightarrow \infty}\eta\phi_{(\mu_{k,\ell},\rho_k)}(x_k,\la_k,s_k)\ge 0$ which implies $\lim_{\ell\rightarrow \infty}\phi_{(\mu_{k,\ell},\rho_k)}(x_k,\la_k,s_k)=0$.

Now we prove the result in case (3). Suppose that $k_i$ and $k_{i+1}$ are the indices of two adjoining iterations such that
\bea \mu_{k_i}>\eta\phi_{(\mu_{k_i},\rho_{k_i})}(x_{k_i},\la_{k_i},s_{k_i}),\quad \mu_{k_{i+1}}>\eta\phi_{(\mu_{k_{i+1}},\rho_{k_{i+1}})}(x_{k_{i+1}},\la_{k_{i+1}},s_{k_{i+1}}),\quad \label{musub}\eea
$\ell_i$ is the number of loops in Step 5.1 of \refal{alg2} such that
\bea   \mu_{k_i,\ell_i}\le\eta\phi_{(\mu_{k_i,\ell_i},\rho_{k_i})}(x_{k_i},\la_{k_i},s_{k_i}).\nn\eea
Since $\mu_{k_i,\ell_i}\ge\gamma\phi_{(\mu_{k_i,\ell_i},\rho_{k_i})}(x_{k_i},\la_{k_i},s_{k_i})$, one has
\bea \mu_{k_i+1}=(1-\alp_{k_i})\mu_{k_i,\ell_i}+\alp_{k_i}\gamma\phi_{(\mu_{k_i,\ell_i},\rho_{k_i})}(x_{k_i},\la_{k_i},s_{k_i})
\le\mu_{k_i,\ell_i}\le\mu_{k_i}/\eta, \nn\eea
and $\mu_{k_{i+1}}\le\mu_{k_i+1}\le\mu_{k_i}/\eta$.
Thus, a strictly monotonically decreasing infinite subsequence $\{\mu_{k_i}\}$ satisfying \reff{musub} is derived. Therefore,
\bea \lim_{i\rightarrow \infty}\mu_{k_i}=0,\quad\lim_{i\rightarrow \infty}\phi_{(\mu_{k_i},\rho_{k_{i}})}(x_{k_{i}},\la_{k_i},s_{k_i})=0. \nn\eea
Note that $\{\mu_k\}$ is a monotonically nonincreasing sequence,
the desired result is straightforward by the preceding equations. \end{proof}

\section{Local convergence}

In this section, we prove that, under suitable conditions, our algorithm with global convergence result (1) of \reft{mainglob} can be quadratically convergent to the KKT point of the original problem. For convenience of statement, we denote $w^*=(x^*,\la^*,s^*)$ and $w_k=(x_k,\la_k,s_k)\in\Re^{2n+m}$ for all $k\ge 0$. The following blanket assumptions are requested for local convergence analysis.

\bas\label{ass2} \ \\
(1) $w_k\rightarrow  w^*$ and $\mu_k\rightarrow  0$ as $k\rightarrow \infty$; \\
(2) The functions $f$ and $h_i\, (i=1,\ldots,m)$ are twice differentiable on $\Re^n$, and their second derivatives are
Lipschitz continuous at some neighborhood of $x^*$;\\
(3) The gradients $\na h_i(x^*) \ (i=1,\ldots,m)$ are linearly independent;\\
(4) There holds $x^*+s^*>0$; \\
(5) $d^TB^*d>0$ for all $d\ne 0$ such that $\na h(x^*)^Td=0$ and $d_j=0$ for $j\in\{j\rvert x_j^*=0, j=1,\ldots,n\}$, where $B^*=\na^2f(x^*)-\sum_{i=1}^m\la_i^{*}\na^2h_i(x^*)$ and $\la^*\in\Re^m$ is the Lagrange multiplier vector associated with  at $x^*$ for all equality constraints, $d_j$ is the $j$-th component of $d$.
\eas

Under \refa{ass2}, $\{s_k\}$ is bounded, thus
$\rho_k$ will keep constant after a finite number of iterations. By \reft{mainglob}, $(x^*,\la^*,s^*)$ is a KKT triple of the original problem.  Without loss of generality,
let $\rho_k=\rho^*$ for all $k\ge 0$, and, correspondingly, $y_k\rightarrow  y^*$ and $z_k\rightarrow  z^*$ as $k\rightarrow \infty$.
It follows from \reff{zydf1} and \reff{zydf} that $z^*=x^*$ and $y^*=s^*/\rho^*$.
Thus, $z_j^*+y_j^*>0$ for all $j=1,\ldots,n$.
\ble\label{NW99} Suppose that \refa{ass2} hold. Let $Y^*=\diag(y^*)$ and $Z^*=\diag(z^*)$. Then the matrix
\bea G^*=\left(\ba{cc}
1+\gamma\frac{D\phi_{(0,\rho^*)}(w^*)}{D\mu} & \gamma(\na_{w}\phi_{(0,\rho^*)}(w^*))^T \\[5pt]
\left(\ba{c}
0 \\
0 \\
-\frac{1}{\rho^*}(Z^*+Y^*)^{-1}e\ea\right) & H^* \ea\right)\nn\eea
is nonsingular, where $\frac{D\phi_{(0,\rho^*)}(w^*)}{D\mu}=\frac{D\phi_{(\mu,\rho)}(w)}{D\mu}\rvert_{(\mu,\rho)=(0,\rho^*),w=w^*}$, \\ $\na_{w}\phi_{(0,\rho^*)}(w^*)=\na_{w}\phi_{(\mu,\rho)}(w)\rvert_{(\mu,\rho)=(0,\rho^*),w=w^*}$, and
\bea H^*=\left(\ba{ccc}
B^* & -\na h(x^*) & -I \\
\na h(x^*)^T & 0 & 0  \\
(Z^*+Y^*)^{-1}Y^* & 0 & \frac{1}{\rho^*}(Z^*+Y^*)^{-1}Z^*\ea\right).\nn\eea
\ele\begin{proof} 
In order to derive the result, we need only to prove that the system
\bea G^*d=0 \nn\eea
has a unique solution $d^*=0$. Corresponding to the partition of $G^*$, $d\in\Re^{2n+m+1}$ has a partition $d=(d_{\mu},d_w)$, where $d_{\mu}\in\Re$, $d_w=(d_x,d_{\la},d_s)$ with $d_x\in\Re^n$, $d_{\la}\in\Re^m$, and $d_s\in\Re^n$. Thus,
\bea   (1+\gamma\frac{D\phi_{(0,\rho^*)}(w^*)}{D\mu})d_{\mu}+\gamma(\na_{w}\phi_{(0,\rho^*)}(w^*))^Td_w=0,  \label{eqn20191212a}\\
  B^*d_x-\na h(x^*)d_{\la}-d_s=0,  \label{eqn20191212b} \\
  \na h(x^*)^Td_x=0,  \label{eqn20191212c}  \\
 -\frac{1}{\rho^*}(Z^*+Y^*)^{-1}ed_{\mu}+(Z^*+Y^*)^{-1}Y^*d_x+\frac{1}{\rho^*}(Z^*+Y^*)^{-1}Z^*d_s=0.  \label{eqn20191212d}\eea
Note that $\frac{D\phi_{(0,\rho^*)}(w^*)}{D\mu}=0$ and $\na_{w}\phi_{(0,\rho^*)}(w^*))=0$ since $z^*-x^*=0$ and $w^*$ is a KKT triple of the original problem.
Thus, due to \reff{eqn20191212a}, $d_{\mu}^*=0$. Furthermore, since $y_j^*z_j^*=0$ for all $j=1,\ldots,n$, \reff{eqn20191212d} implies $(d_x^*)^Td_s^*=0$, and $(d_x^*)_j=0$ when $x_j^*=0$, $(d_s^*)_j=0$ as $s_j^*=0$, where $(d_x^*)_j$ and $(d_s^*)_j$ are, respectively, the $j$-th components of $d_x^*$ and $d_s^*$. Hence, \bea
(d_x^*)^TB^*d_x^*=0,\quad \na h(x^*)^Td_x^*=0,\quad (d_{x}^*)_j=0\ \hbox{for}\ j\in\{j\rvert x_j^*=0, j=1,\ldots,n\}, \nn \eea
which, due to \refa{ass2} (5), implies $d_x^*=0$. Finally, $d_{\la}^*=0$ follows from \refa{ass2} (3) since $\na h(x^*)d_{\la}=0$.
\end{proof} 

The preceding proof also shows that $H^*p=0$ implies $p=0$. Thus, $H^*$ is also nonsingular. Let $w=(x,\lambda,s)$ and $$\Phi(\mu,w)=\left(\ba{c}
\mu+\gamma\phi_{(\mu,\rho^*)}(w) \\[5pt]
\na f(x)-\na h(x)\lambda-s\\
h(x)\\
z-x\ea\right).$$
Then $\Phi(0,w^*)=0$.
The following lemma can be obtained in a way similar to Lemma 2.1 in \cite{ByrLiN97}.
We will not give its proof for brevity.
\ble\label{lemc1} Suppose that Assumption \ref{ass2} holds. Then there are sufficiently small scalar $\epsilon>0$
and positive constants $M_0$ and $L_0$, such that,
for all $(\mu,w)\in\{(\mu,w)\in\Re_{++}\times\Re^{2n+m}\rvert\|(\mu,w)-(0,w^*)\|<\epsilon\}$, $\na_{(\mu,w)} \Phi(\mu,w)$ is invertible,
$\|[\na_{(\mu,w)}  \Phi(\mu,w)]^{-1}\|\le M_0,$
and
\bea\|(\na_{(\mu,w)} \Phi(\mu,w))^T((\mu,w)-(0,w^*))-\Phi(\mu,w)\|\le L_0\|(\mu,w)-(0,w^*)\|^2,\quad \label{20140415d}\eea 
where $\na_{(\mu,w)} \Phi(0,w^*)=\na_{(\mu,w)} \Phi(\mu,w)\rvert_{\mu=0,w=w^*}$.\ele

Using \refl{lemc1},
the following result shows that the step $(\De\mu_k, d_k)$ can be a quadratically or superlinearly convergent step.
\begin{thm} \label{localt} Suppose that Assumption \ref{ass2} holds. Then there is a sufficiently small scalar $\epsilon>0$, such that,
for all $(\mu_k,w_k)\in\{(\mu,w)\in\Re_{++}\times\Re^{2n+m}\rvert\|(\mu,w)-(0,w^*)\|<\epsilon\}$, one has the following results.\\
(1) If $\|(B_k-B^*)d_x\|$ $=O(\|d_x\|^2)$ for every $d_x\in\Re^n$, then
\bea {\|(\mu_k,w_{k})+(\Delta\mu_k,d_{k})-(0,w^*)\|}=O{(\|(\mu_k,w_k)-(0,w^*)\|^2)}. \eea
That is, $(\Delta\mu_k, d_{k})$ is a quadratically convergent step. \\
(2) If $\|(B_k-B^*)d_x\|$ $=o(\|d_x\|)$ for every $d_x\in\Re^n$, then
\bea \|(\mu_k,w_{k})+(\Delta\mu_k,d_{k})-(0,w^*)\|=o{(\|(\mu_k,w_k)-(0,w^*)\|)}, \eea
i.e., $(\Delta\mu_k, d_{k})$ is a superlinearly convergent step. \end{thm} \begin{proof} 
In order to prove the result (1), we show
\bea\limsup_{k\rightarrow \infty}{\|(\mu_k,w_k)+(\Delta\mu_k,d_{k})-(0,w^*)\|}/{\|(\mu_k,w_k)-(0,w^*)\|^2}\le\xi,
\label{local1}\eea
where $\xi>0$ is a constant.

Let $\Phi_k=\Phi(\mu_k,w_k)$, $J_k=\na_{(\mu,w)}\Phi(\mu_k,w_k)^T$, $G_k$ is a matrix which has the same components as $J_k$ except that the Lagrangian Hessian $\na^2_{xx}L(w_k)=\na^2 f(x_k)-\sum_{i=1}^m\la_{ki}\na^2 h_i(x_k)$ in $J_k$ is replaced by $B_k$.
Then $G_k(\Delta\mu_k,d_{k})=-\Phi_k$. By \refl{lemc1}, $J_k$ is invertible. Note that \bea
G_k=J_k+G_k-J_k=J_k+\left(\ba{ccc}
0 & 0 & 0 \\
0 & B_k-B^* & 0 \\
0 & 0 & 0 \ea\right)-\left(\ba{ccc}
0 & 0 & 0 \\
0 & \na^2_{xx}L(w_k)-B^* & 0 \\
0 & 0 & 0 \ea\right), \nn\eea
it follows from the condition $\|(B_k-B^*)d_x\|$ $=O(\|d_x\|^2)$ and \refa{ass2} (2) that $G_k$ is invertible and $\|G_k^{-1}\|\le M_0$ for some scalar $M_0>0$ and for all sufficiently large $k$.
Thus, $\|(\Delta\mu_k,d_{k})\|=O(\|(\mu_k,w_k)-(0,w^*)\|)$. Moreover, \bea
G_k(\Delta\mu_k,d_{k})=J_k(\Delta\mu_k,d_{k})+(B_k-B^*)d_{xk}-(\na^2_{xx}L(w_k)-B^*)d_{xk}=-\Phi_k. \nn\eea
Therefore,
\bea   
&&\|(\mu_k,w_k)+(\Delta\mu_k,d_{k})-(0,w^*)\|\nn\\
  &&= \|J_k^{-1}(J_k((\mu_k,w_k)-(0,w^*))-\Phi_k-(B_k-B^*)d_{xk}+(\na^2_{xx}L(w_k)-B^*)d_{xk})\|\nn\\
  &&\le M_0[L_0\|(\mu_k,w_k)-(0,w^*)\|^2+O(\|(\mu_k,w_k)-(0,w^*)\|^2)], \label{20140415g}
\eea
where the last inequality follows from \reff{20140415d} of \refl{lemc1}.
Thus, \reff{local1} follows immediately from \reff{20140415g}.

If $\|(B_k-B^*)d_x\|$ $=o(\|d_x\|)$, then the last inequality \reff{20140415g} should be \bea
  &&\|(\mu_k,w_k)+(\Delta\mu_k,d_{k})-(0,w^*)\|\nn\\
  &&\le M_0[L_0\|(\mu_k,w_k)-(0,w^*)\|^2+o(\|(\mu_k,w_k)-(0,w^*)\|)]. \label{20200102a}
\eea
Hence, the result (2) follows immediately.
\end{proof} 

Now we prove that, under suitable conditions, our algorithm can be quadratically convergent to the KKT triple of the original problem.
\begin{thm}  Suppose that Assumption \ref{ass2} holds.
If $\|(B_k-B^*)d_x\|$ $=O(\|d_x\|^2)$ for every $d_x\in\Re^n$,
$\tau<1/2$, then either $\mu_{k+1}=\mu_k$ or $\mu_{k+1}=\gamma_0\phi_{(\mu_k,\rho^*)}(w_k)$, $x_{k+1}=x_k+d_{xk}$,
$s_{k+1}=s_k+d_{sk}$, and $\la_{k+1}=\la_k+d_{\la k}$ for all sufficiently large $k$.
Moreover, ${\|w_{k+1}-w^*\|}=O{(\|w_k-w^*\|^2)}$.
\end{thm} 
We need to prove that, for all sufficiently large $k$, $\alp_k=1$ will be accepted by the line search procedure \reff{steprulea}. By \reft{localt}, 
\bea
  &&\phi_{(\mu_k+\De\mu_k,\rho^*)}(w_k+d_k) \nn\\
  &&=\phi_{(\mu_k+\De\mu_k,\rho^*)}(w_k+d_k)-\phi_{(0,\rho^*)}(w^*) \nn\\
  &&=(\na_{(\mu,w)}\phi_{(\mu_k+\De\mu_k,\rho^*)}(w_k+d_k))^T((\mu_k,w_k)+(\De\mu_k,d_k)-(0,w^*))\nn\\
  &&\quad+O(\|(\mu_k,w_k)+(\De\mu_k,d_k)-(0,w^*)\|^2) \nn\\
  &&=O(\|(\mu_k,w_k)-(0,w^*)\|^2). \nn
  \eea
Note that $\phi_{(\mu_k,\rho^*)}(w_k)=\phi_{(\mu_k,\rho^*)}(w_k)-\phi_{(0,\rho^*)}(w^*)=O(\|(\mu_k,w_k)-(0,w^*)\|)$ and $\tau<\hf$. Thus, 
\bea (1-2\tau)\phi_{(\mu_k,\rho^*)}(w_k)=O(\|(\mu_k,w_k)-(0,w^*)\|), \nn\eea
and the full step will be accepted by \reff{steprulea}.

By \reft{localt} (1), \bea {\|(\mu_{k+1},w_{k+1})-(0,w^*)\|}=O{(\|(\mu_k,w_k)-(0,w^*)\|^2)}. \label{20200407a}\eea
Due to $\mu_k\le\eta\phi_{(\mu_k,\rho^*)}(w_k)$, $\mu_k=O(\|w_k-w^*\|)$. The desired result follows from \reff{20200407a} immediately.

\section{Numerical experiments}

Our method can be easily extended to solve the nonlinear programs with general equality and inequality constraints
\begin{align} 
\hbox{min}   &~f(x) \label{nlp1-1}\\
\st   &~h(x)=0,\quad g(x)\ge 0, \label{nlp1-2} 
\end{align} 
by substituting \reff{zydf1} and \reff{zydf}, respectively, with \bea
   z_j(x,s;\mu,\rho)\equiv\frac{1}{2\rho}\Big(\sqrt{(s_j-\rho g_j(x))^2+4\rho\mu}-(s_j-\rho g_j(x))\Big), \nn\\
   y_j(x,s;\mu,\rho)\equiv\frac{1}{2\rho}\Big(\sqrt{(s_j-\rho g_j(x))^2+4\rho\mu}+(s_j-\rho g_j(x))\Big), \nn\eea
where $g: \Re^n\rightarrow \Re^{m_{\cal I}}$ is a twice continuously differentiable real-valued function on $\Re^n$, $j=1,\ldots,m_{\cal I}$.
No slack variables are introduced to cope with the general inequality constraints, which is different from the technique commonly used in interior-point methods for nonlinear programs \reff{nlp1-1}--\reff{nlp1-2}.
Our numerical experiments are conducted on a Lenovo laptop with the LINUX operating system (Fedora 11).
\refal{alg2} is implemented in MATLAB (version R2008a).

The algorithm is firstly used to solve a well-posed nonlinear program from the literature. The test problem was presented by W\"{a}chter and Biegler \cite{WacBie00} and further discussed by Byrd, Marazzi and Nocedal \cite{ByrMaN01}:
\begin{eqnarray}
\hbox{min} && x \label{test1obj}\\
\hbox{s.t.} && x\sp{2}-1\ge 0, \quad x-2\ge 0. \label{test1cons}
\end{eqnarray}
This problem is well-posed since it has a unique global minimizer $x^*=2$, at which both the linear independence constraint qualification (short for LICQ) and the Mangasarian-Fromowitz constraint qualification (short for MFCQ) hold. However, starting from $x_0=-4$, \cite{WacBie00} showed that many line-search infeasible interior-point methods may be jammed and fail to find the solution.

\refal{alg2} is then used to find the solutions for a set of nonlinear programming test problems of the CUTEr collection \cite{BonCGT95}. Since the code is very elementary, we restricted our test problems to the $122$ HS problems, where problems HS101--103 were excluded since they are only defined on positive variables.
These test problems include not only the problems with general equality and inequality constraints, but also the problems with bound constraints and the problems with only equality constraints \cite{HocSch81}.

In our implementation, the initial parameters are selected as follows: $\mu_0=0.1$, $\rho_0=1$, $\eta=10$, $\gamma_0=0.001$, $\delta=0.5$, $\tau=0.01$, $\sigma=0.01$, and $\epsilon=10^{-8}$. For all $k\ge 0$, we take $B_k$ to be the exact Lagrangian Hessian provided that it is positive semi-definite (where the gradient and Hessian are provided by the test sets). Otherwise, we modify $B_k$ to $B_k+\xi I$ with $\xi>0$ being as small as possible so that the modified Hessian is positive semi-definite.

For comparison, these test problems are also solved by the well regarded and recognized interior-point solver IPOPT \cite{WacBie06} (Version 3.0.0). In implementation, \refal{alg2} can use the KKT residuals of the original problem directly as the measure of our terminating conditions: \bea
E(x_k,\lambda_k,s_k)\le\epsilon, \label{20200413}\eea where 
\begin{align*}
E(x_k,\lambda_k,s_k)=
\max\{&\|\na f(x_k)-\na h(x_k)\lambda_k-s_k\|_{\infty},\|h(x_k)\|_{\infty},\\
&\|\max\{-(x_k+s_k),0\}\|_{\infty},\|x_k\circ s_k\|_{\infty}\},
\end{align*}
$x_k\circ s_k$ is the Hadamard product of $x_k$ and $s_k$. If one has the scaling parameters $s_d=1$ and $s_c=1$ in the terminating conditions of \cite{WacBie06},
then the accuracy differences between \refal{alg2} and IPOPT should be in the range of the tolerance.

For test problem \reff{test1obj}--\reff{test1cons}, we use the standard initial point $x_0$ as the starting point,
and set $s_0$ to be the all-one vector. The implementation of our algorithm terminates at  
$x^*=2$ together with $s_1^*=-1.1972\times 10^{-16}$, $s_2^*=1.0000$ in $4$ iterations. Both the numbers of function and gradient evaluations are $5$.
See Table \ref{tab1} for more details on iterations. From there one can observe the rapid convergence of $\mu_k$,
$\phi_{(\mu_k,\rho_k)}(x_k,\lambda_k,s_k)$ and $E(x_k,\lambda_k,s_k)$, where $\mu_k$ is the current value of the parameter, $x_k$ and $s_k$ are the estimates of the primal and dual variables, respectively,
$f_k=f(x_k)$, $v_k$ is the $\ell_{\infty}$ norm of violations of constraints, $\phi_k=\phi_{(\mu_k,\rho_k)}(x_k,\lambda_k,s_k)$,
$E_k=E(x_k,\lambda_k,s_k)$. As a comparison, IPOPT fails to find the solution and terminates at $x^*=-1.0000$ in $13$ iterations. In interior-point framework, this problem has been solved by the recently developed methods
of \cite{DLS17} and \cite{LiuDai18} in totally 16 and 19 iterations, respectively.
\begin{table}[ht]
	\caption{Output of \refal{alg2} for test problem \reff{test1obj}--\reff{test1cons}}\label{tab1} 
	\renewcommand\arraystretch{1.5}
	\setlength{\tabcolsep}{1.ex}
	\noindent\[
	\begin{tabular}{|c|c|c|c|c|c|c|c|}
	\hline
	$k$ & $\mu_k$ & $x_k$ & $s_k$ & $f_k$ & $v_k$ & $\phi_k$ & $E_k$ \\
	\hline
	\hline
	0 & 0.1 & -4 & $(1,1)$ & -4 & 6 & 50.5785 & 15 \\
	1 & 0.0506 & 2.0190 & $(0.0276,1.2212)$ & 2.0190 & 0 & 0.0557 & 0.3328 \\
	2 & 5.5681e-05 & 2.0080 & $(0.0002,0.9992)$ & 2.0080 & 0 & 3.1754e-05 & 0.0080 \\
	3 & 3.1754e-13 & 2.0000 & $(0.0000,1.0000)$ & 2.0000 & 4.6437e-07 & 1.2875e-13 & 6.1372e-07 \\
	4 & 1.2875e-16 & 2 & $(-0.0000,1.0000)$ & 2 & 0 & 6.1630e-33 & 3.5916e-16 \\
	\hline
	\end{tabular}
	\]
\end{table}

When solving the HS test problems of the CUTEr collection, \refal{alg2} was terminated as either $E(x_k,\lambda_k,s_k)\le\epsilon$, or the number of iterations is larger than 500 (which is the default setting of IPOPT), the step-size is too small ($\alp_k\le\delta^{40}$), the coefficient matrix of the system \reff{dsysa} is degenerate. The latter three cases of termination can be resulted from that the Hessian does not satisfy \refa{ass1} (3), the condition (4) of \refa{ass1} does not hold, and some test problems are only defined on strictly positive variables.

Since we do not require the iterates to be interior points, our algorithm has the freedom to use the standard initial points for all HS test problems. However, for the purpose of comparison, we have modified the initial points in line with the initialization of IPOPT \cite{WacBie06}. In our implementation, \refal{alg2} successfully solved $79$ problems and terminated with \reff{20200413}, while IPOPT found the approximate solutions of $121$ problems satisfying its default terminating conditions, where only for problem HS87 IPOPT reached its restriction of the maximum of the total number of iterations.

In order to further observe how \refal{alg2} performs in solving nonlinear programming test problems, we provide $4$ figures Figures \ref{fig1}--\ref{fig4} to show log scaling performance profiles (see Dolan and Mo\'re \cite{DolMor}) of our algorithm in comparison with IPOPT on both solved $79$ problems with respect to iteration count, function evaluations, gradient evaluations, and the CPU time, where IPRM represents our primal-dual interior-point relaxation method (\refal{alg2}), respectively. Figures \ref{fig1}--\ref{fig3}  show that, under the measures on the former three items, IPRM performs approximate but inferior to IPOPT. However, Figure \ref{fig4}  shows that IPRM needs less CPU time than IPOPT, which may be partially resulted from that the system \reff{dsysa} in IPRM is solved by the MATLAB's built-in ``backslash" command and that our algorithm does not incorporate any sophisticated techniques such as inertia correction, feasibility restoration, and so on.

Since our method is currently at a very early stage of development, and we note that a nonmonotone line search variant of our algorithm can successfully solve more than $100$ HS test problems, 
it is not surprising that our implementation of \refal{alg2} is not very convincing in comparison to the very regarded and recognized IPOPT. However, it is still encouraging by the numerical experiments since \refal{alg2} has still much space for improvement such as incorporating some scaling and inertial control techniques and using some robust subroutine and solver for solving the system \reff{dsysa} more efficiently.

\bfig[ht!b]
\centering
\includegraphics[width=0.75\textwidth,height=0.5\textwidth]{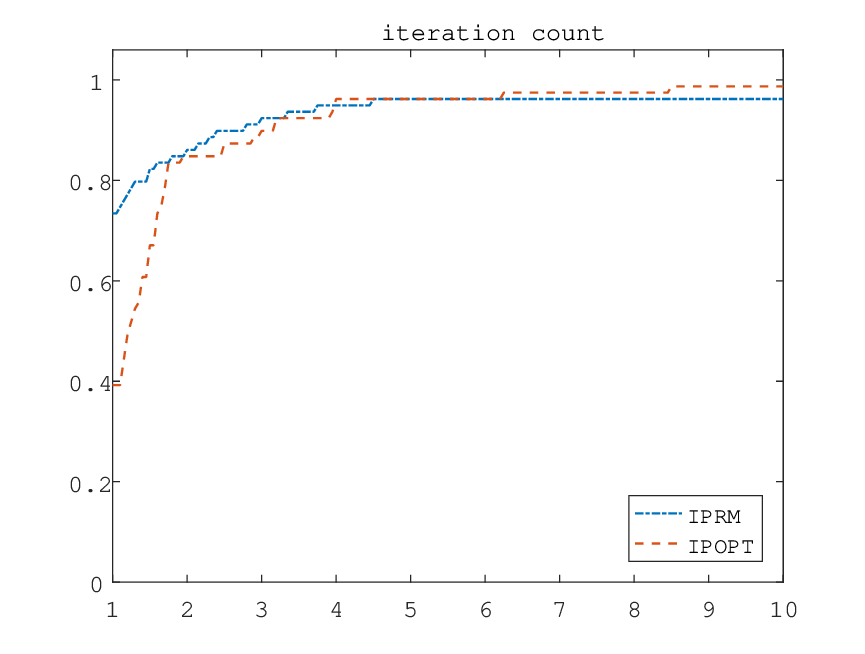}\\
\caption{Performance plot for iteration count}\label{fig1}
\efig

\begin{figure}[ht!b] 
	\centering
	\includegraphics[width=0.75\textwidth,height=0.5\textwidth]{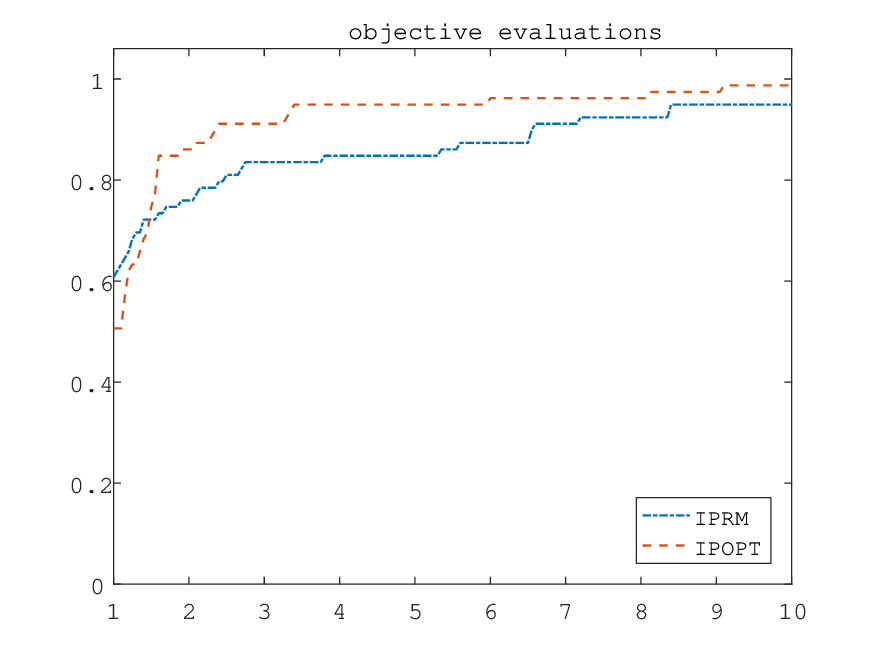}\\
	\caption{Performance plot for function evaluations}\label{fig2}
\end{figure}

\begin{figure}[ht!b] 
	\centering
	\includegraphics[width=0.75\textwidth,height=0.5\textwidth]{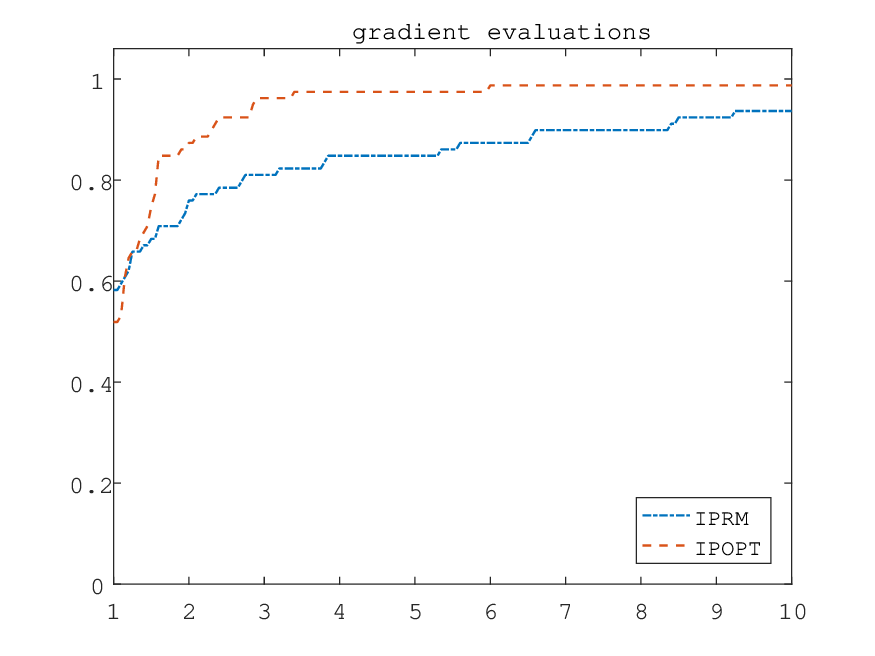}\\
	\caption{Performance plot for gradient evaluations}\label{fig3}
\end{figure}

\begin{figure}[ht!b] 
	\centering
	\includegraphics[width=0.75\textwidth,height=0.5\textwidth]{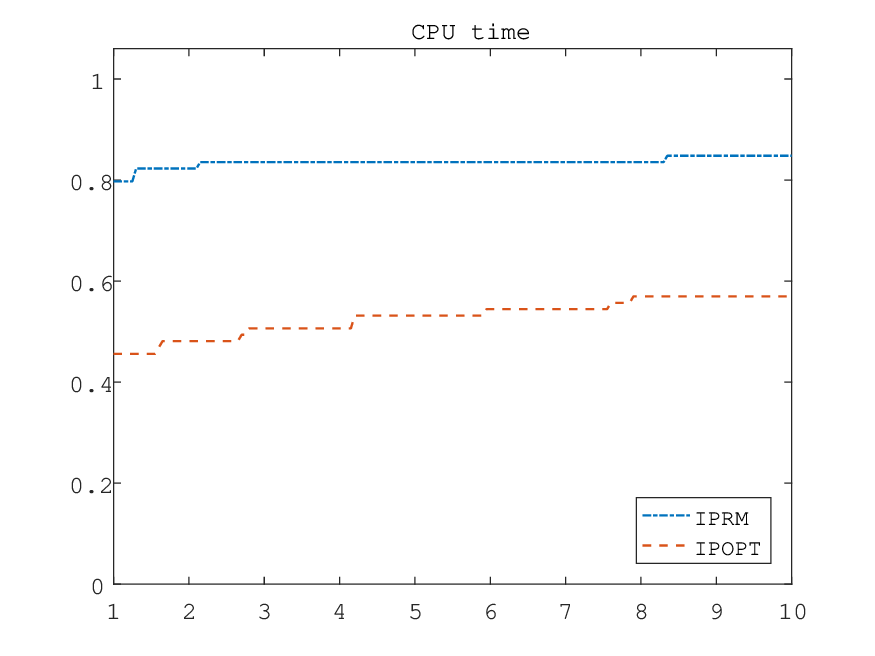}\\
	\caption{Performance plot for the CPU time}\label{fig4}
\end{figure}


\section{Conclusion}

We present a novel primal-dual interior-point relaxation method for nonlinear programs with general equality and nonnegative constraints in this paper. The method can be easily extended to solve the problems with general inequality constraints. It is based on solving a parametric equality constrained mini-max subproblem. Our method is of the interior-point variety, but does not require any primal or dual iterates to be interior. A new smoothing approach is introduced. Our method is capable of circumventing the jamming difficulty which results in that many interior-point methods failed to converge to the solution and improving the ill conditioning of the classic primal-dual interior-point methods as the barrier is small. Under suitable conditions, our method is proved to be globally convergent and locally quadratically convergent to the KKT triple of the original problem. Preliminary numerical results on a well-posed problem for which many line-search interior-point methods fail to find the minimizer and a set of test problems from CUTEr collection show that our method is efficient.

\begin{acknowledgements}
The research is supported by the NSFC grants (nos. 12071108, 11671116, 12021001,
11991021, 11991020, 11971372, and 11701137), National Key R\&D Program of China (nos.
2021YFA1000300 and 2021YFA1000301), the Strategic Priority Research Program of Chinese Academy of
Sciences (no. XDA27000000), and the Natural Science Foundation of Hebei Province (no. A2021202010).
\end{acknowledgements}

\end{document}